\documentclass[10pt]{amsart}
\usepackage{amsfonts}
\usepackage{amsmath}
\usepackage{amssymb}
\usepackage{amsmath,amssymb,amsfonts,amsthm,graphics,
latexsym, amscd, amsfonts, epsfig, eepic,epic}
\usepackage{mathrsfs}
\usepackage{color}
\usepackage{eucal}
\usepackage{eufrak}
\usepackage[all]{xypic}
\usepackage{xspace}

\newtheorem{theorem}{Theorem}[section]

\newtheorem{proposition}[theorem]{Proposition}
\newtheorem{corollary}[theorem]{Corollary}

\theoremstyle{remark}

\makeatletter \@addtoreset{equation}{section}

\title[RIP: partition function and base pair
       pairing probabilities]
      {RNA-RNA interaction prediction:\\ partition function and base pair
       pairing probabilities}
\author{Fenix W.D. Huang$^1$, Jing Qin$^1$, Christian M.
Reidys$^{1,\star}$ and Peter F. Stadler$^2$}
\address{$^1$Center for Combinatorics, LPMC-TJKLC 
          \\
         Nankai University  \\
         Tianjin 300071\\
         P.R.~China\\
         Phone: *86-22-2350-6800\\
         Fax:   *86-22-2350-9272\\
         $^2$University of Leipzig \\
         Leipzig D-04107\\
         H\"artelstr. 16-18\\
         Germany\\}
\email{duck@santafe.edu}
\thanks{}
\date{March, 2009}
\begin{document}
\maketitle
\begin{abstract}
In this paper, we study the interaction of an antisense RNA and its
target mRNA, based on the model introduced by Alkan {\it et al.}
(Alkan {\it et al.}, J.~ Comput.~ Biol., Vol:267--282, 2006). Our
main results are the derivation of the partition function
 \cite{Backhofen} (Chitsaz {\it et al.}, Bioinformatics, to appear,
2009), based on the concept of tight-structure and the computation
of the base pairing probabilities. This paper contains the folding
algorithm {\sf rip} which computes the partition function as well as
the base pairing probabilities in $O(N^4M^2)+O(N^2M^4)$ time and
$O(N^2M^2)$ space, where $N,M$ denote the lengths of the interacting
sequences.
\end{abstract}

\maketitle

\section{Introduction}\label{S:Introduction}

The discovery of small RNAs that bind to their target mRNAs in order
to prohibit their translation and down-regulate the expression
levels of corresponding genes has drawn a lot of attention in the
RNA world \cite{McManus}. Studies have shown that many RNA-RNA
interactions play a significant role in different cellular
processes, such as mediate pseudouridylation and methylation of rRNA
\cite{Bachellerie}, nucleotide insertion into mRNAs \cite{Benne},
splicing of pre-mRNA \cite{Zorio} and translation control or plasmid
replication control \cite{Banerjee, Fire, Kugel}.

Regulatory RNAs constitute a subclass of the antisense RNA family;
encompassing the snRNAs, gRNAs and snoRNAs that play a role in the
context of rRNA modification, RNA editing, mRNA spicing and plasmid
copy-number regulation. In addition, antisense RNAs are synthesized
for studying specific gene functions. Since the first published
result on natural antisense RNAs which regulate gene expression in
C.~elegans \cite{Parrish,Yang,Hammond,Reinhart}, Drosophila
\cite{Nykanen}, and other organisms \cite{Wagner}, the problem of
predicting how two nucleic acid strands interact--the so called
RNA-RNA interaction problem (RIP)--has come into focus.

As observed by Alkan {\it et al.} \cite{Alkan:06}, the RIP is
NP-complete. The actual argument constitutes an extension of the
work of Akutsu \cite{Akutsu} derived in the context of single RNA
secondary structure prediction problems with pseudoknots. As in
Rivas and Eddys pseudoknot folding algorithm \cite{Rivas} the
general idea here is to consider specific classes of interactions,
that can be computed via dynamic programming routines.
There are several other methods that consider somewhat restricted
versions of the RNA-RNA interaction. For instance, one method
concatenates the two interacting sequences and subsequently employs
a slightly modified standard secondary structure folding algorithm.
The algorithms RNAcofold \cite{Hofacker, Bernhart}, pairfold
\cite{Andronescu} and NUPACK \cite{Ren} subscribe to this strategy.
However, this approach cannot predict important motifs in RIPs, as
for instance kissing hairpin loops. The concatenation idea has also
been employed using the pseudoknot folding algorithm of Rivas and
Eddy \cite{Rivas}. The resulting algorithm, however, does still not
generate all relevant interaction structures
\cite{Backhofen,Reidys:frame}. An alternative line of thought is to
neglect all internal base-pairings in either strand and to compute
the minimum free energy (mfe) secondary structure for their
hybridization under this constraint. For instance, RNAduplex follows
this line of thought making it formally equivalent to the classic
secondary structure folding algorithm of Waterman
\cite{Waterman:78,Waterman:80,Waterman:86,Waterman:94a}. Furthermore
we have the algorithm RNAup \cite{Mueckstein:05a, Mueckstein:08a}
which uses the Alkan's model, allowing for {\it one} interaction
region having unbranched interactions within any loop. RNAup can
therefore capture single but not multiple kissing hairpins. Finally
there is IntaRNA \cite{Busch:08} facilitating the efficient
prediction of bacterial sRNA targets incorporating target site
accessibility and seed regions.

\begin{figure}[ht]
\centerline{%
\epsfig{file=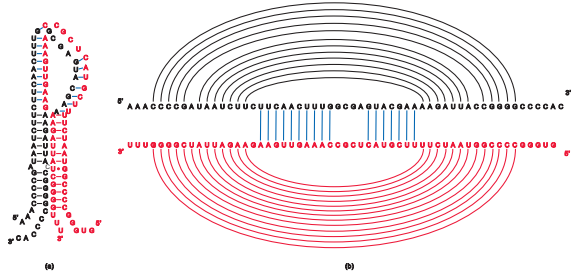,width=0.5\textwidth}\hskip15pt }
\caption{\small Natural joint structure between small RNA molecules
CopA(antisense) and CopT(target) in E.coli \cite{Alkan:06}.}
\label{F:interaction}
\end{figure}

Alkan {\it et al.} \cite{Alkan:06} derived a mfe algorithm for
predicting the joint secondary structure of two interacting RNA
molecules with polynomial time complexity. Here ``joint structure'',
see Fig.~\ref{F:interaction} for example, means that the
intramolecular structures of each molecule are pseudoknot-free, the
intermolecular binding pairs are noncrossing and there exist no so
called ``zig-zags'' (see Section~\ref{S:Introduction} for details).
Zig-zags are sometimes referred to as tangles.

\begin{figure}[ht]
\centerline{%
\epsfig{file=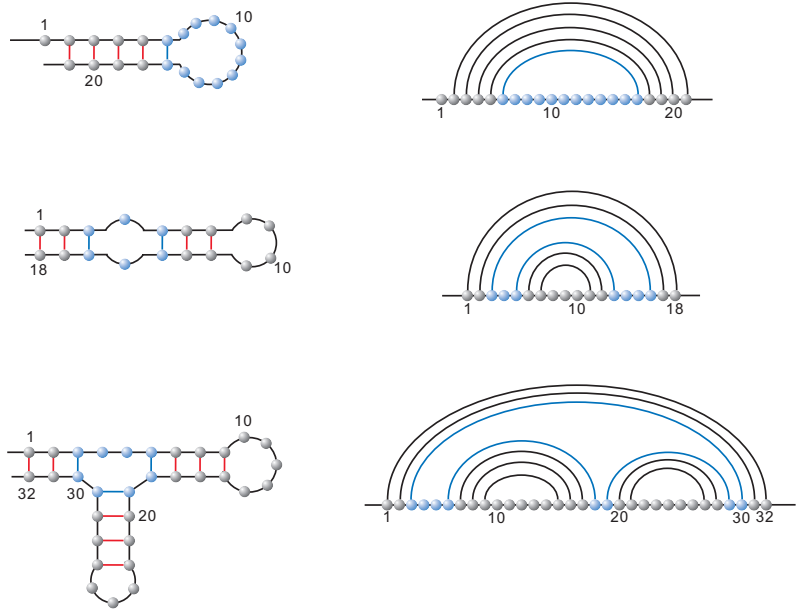, width=0.7\textwidth}\hskip15pt
 }
\caption{\small The standard loop-types for RNA secondary structures:
hairpin-loop (top),
interior-loop (middle) and multi-loop (bottom).} \label{F:stand}
\end{figure}

\begin{figure}[ht]
\centerline{%
\epsfig{file=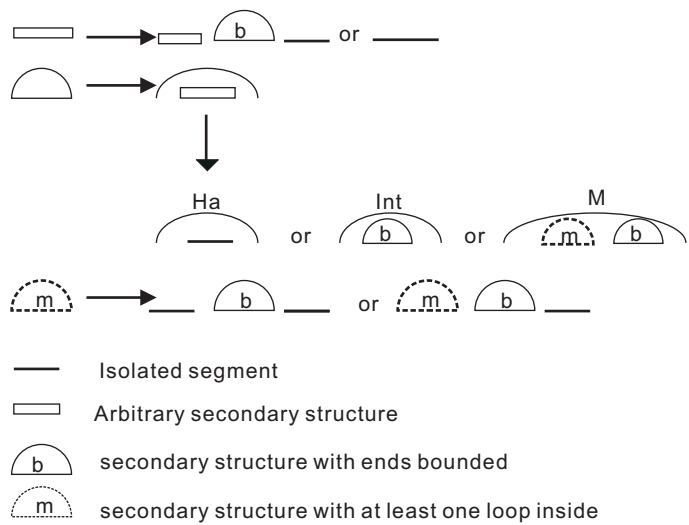,width=0.6 \textwidth}\hskip15pt
 }
\caption{\small The unique decomposition of secondary structures.}
\label{F:parsec}
\end{figure}

Recently, Chitsaz et.al. \cite{Backhofen} presented a dynamic
programming algorithm which computes the partition function in
$O(N^6)$ time. The key point for passing from the mfe folding of
Alkan \cite{Alkan:06} to the partition function is a unique grammar
by which each interaction structure can be generated. The dynamic
programming routine for the partition function of RNA secondary
structures is due to McCaskill \cite{McCaskill} and can be outlined
as follows: the free energy of a secondary structure is assumed
additive in terms of its loops $F(S)=\sum_{L\in S}F_{L}$, where
$F_{L}$ denotes the free energy of a loop,$L$. The additivity of the
free energy translates itself into the multiplicativity in the
contributions to the partition function $Q$ defined by
$Q=\sum_{S}e^{-F(S)/kT}$, where $Q$ is the sum over all the
secondary structures $S$ of length $M$. This factorization of terms
can be realized by introducing $Q^{b}(i,j)$, where the sum is taken
over all substructures $S[i,j]$ on the segment $[i,j]$ for which
$(S[i],S[j])\in S[i,j]$ and $Q^{s}(i,j)$ for all the configurations
on $[i,j]$, irrespective of whether or not $i,j$ are connected. In
particular, we have $Q^{s}(1,M)=Q$. Consequently, we arrive at the
recursion, see Fig.~\ref{F:parsec}
\begin{equation}
Q^{s}(i,j)=1+\sum_{h,\ell}Q^{s}(i,h-1)Q^{b}(h,\ell).
\end{equation}
Let us next recall the basic loops-types upon which the partition
function and energy parameters \cite{Mathews} of RNA secondary
structures are based:\\
{\bf (1)} a {\it hairpin}-loop (${\sf Ha}(i,j)$), is a pair
$((i,j),[i+1,j-1])$, where $(i,j)$ is an arc and $[i+1,j-1]$ is an
interval, i.e.~a sequence of consecutive vertices
$(i,i+1,\dots,j-1,j)$, having
energy parameter $e^{-G^{\sf Ha}(i,j)/kT}$. \\
{\bf (2)} an {\it interior}-loop (${\sf Int}(i_1,j_1;i_2,j_2)$), is
a sequence $((i_1,j_1),[i_1+1,i_2-1],(i_2,j_2),[j_2+1,j_1-1])$,
where $(i_2,j_2)$ is nested in $(i_1,j_1)$ having the energy parameter
$e^{-G^{\sf Int}(i_{1},j_{1};i_{2},j_{2})/kT}$\\
{\bf (3)} a {\it multi}-loop (${\sf M}(i_0,j_0)$), see
Fig.\ref{F:stand}, is a sequence
\begin{equation}
([i_{0},i_1-1],((i_1,j_1),[i_1+1,j_1-1]),\dots,
((i_t,j_t),[i_t+1,j_t-1]),[j_t+1,j_0])
\end{equation}
having energy parameter
$e^{-(\alpha_{1}+\alpha_{2}(t+1)+\alpha_{3}c_2)/kT}$, where
$\alpha_1, \alpha_2, \alpha_3\in \mathbb{R}$, $t$ is the number of
$R[i_0+1,j_0-1]$-maximal arcs inside $R[i_{0},j_{0}]$ and $c_2$ is
the
number of isolated vertices contained in $[i_0,j_0]$.\\
Based on the above loop-energies, we obtain the following recursion
for $Q^{b}(i,j)$
\begin{align*}
Q^b(i,j)&=e^{-G^{\sf{Ha}}(i,j)/kT}+
\sum_{k_1,k_2}e^{-G^{\sf{Int}}(i,j,k_1,k_2)/kT}\\
&\qquad \qquad\qquad\; \quad +\sum_{\ell}Q^{m1}(i+1,\ell)Q^{
m}(\ell+1,j-1)e^{-(\alpha_1+2\alpha_2)/kT},
\end{align*}
where
\begin{eqnarray*}
Q^{m1}_{i,j} & = &
\sum_{i\le \ell<j}Q^b(k,j)e^{-(\alpha_2+\alpha_3(\ell-i))/kT}\\
Q^{m}(i,j) & = & \sum_{i\le \ell<j}Q^{m1}(i,k)(Q^{m}(\ell+1,j)+
e^{-\alpha_3(j-\ell)/kT}).
\end{eqnarray*}

The key idea in this paper, which eventually leads to the derivation
of both: the partition function as well as the base pairing
probabilities, is the concept of a ``tight structure'', introduced in
Section~\ref{S:comb}. The tight structure plays a central role in
our grammar and is the main tool for obtaining the base pairing
probabilities. This paper includes the folding algorithm {\sf
rip}, which derives the partition function as well as the base
pairing probabilities in $O(N^4M^2)+O(N^2M^4)$ time and $O(N^2M^2)$
space. The source code of {\sf rip} is available upon request.


\section{Combinatorics of interaction structures}\label{S:comb}
In this section we discuss some combinatorial properties of RNA
interaction structures. The key idea introduced here is that of
a tight structure. The main results of this section are: \\
$\bullet$ there exist only four ``types'' of tight structures \\
$\bullet$ given a joint structure $J(i,j;h,\ell)$, each interaction
          bond $(R[i_{0}],S[j_{0}])\in J(i,j;h,\ell)$ is
          contained in a unique $J(i,j;h,\ell)$-tight structure \\
$\bullet$ each joint structure uniquely decomposes into a sequence
of tight structures and secondary structure segments\\
$\bullet$ there exists a unique (but not canonical) decomposion of
a tight structure.\\

\begin{figure}[ht]
\centerline{%
\epsfig{file=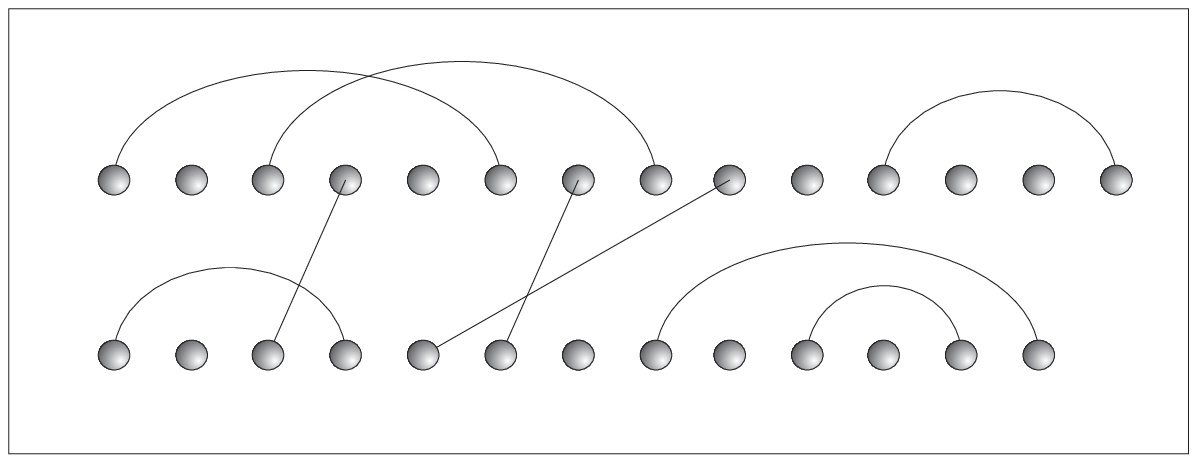,width=0.5\textwidth}\hskip15pt }
\caption{\small A complex $C$ induced by $R[1,14]$ and $S[1,13]$.}
\label{F:complex}
\end{figure}

Let us begin by making precise what we mean by interaction
structures. Suppose we are given two diagrams
\cite{Reidys:07pseu,Reidys:07lego,Reidys:08tan,Reidys:09eff}, $R$
and $S$ of length $N$ and $M$, respectively. Let $R[i]$ and $S[i]$
denote the vertex $i$ of $R$ and $S$, respectively. We shall assume
that $R[1]$ denotes the $5'$ end of $R$ and $S[1]$ denotes the $3'$
end of $S$ as RNA sequences. The induced subgraph of $S$ with
respect to the subsequence $(S[i],\dots,S[j])$ is denoted by
$S[i,j]$. In particular, $S[i,i]=S[i]$ and $S[i,i-1]=\varnothing$. A
complex $C(R,S,I)$ is a graph consisting of $R,S$ and a set of arcs
of the form $(R[i],S[j])$, $I$, see Fig. \ref{F:complex}. We shall
represent a complex $C$ by drawing $R$ on top of $S$ with the
$R$-arcs in the upper, the $S$-arcs in the lower halfplane and
$I$-arcs vertical. Given a complex $C$, a subcomplex is the subgraph
of $C$, induced by $R[i_1,j_1]$ and $S[i_2,j_2]$.

\begin{figure}[ht]
\centerline{%
\epsfig{file=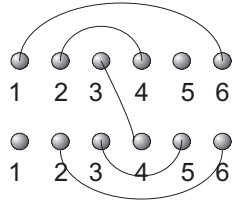,width=0.2\textwidth}\hskip15pt
 }
\caption{\small Ancestors and parents:
for the exterior arc $(R[3],S[4])$, we have the following ancestor sets
$A_{R}(R[3],S[4])=\{(R[1],R[6]),(R[2],R[4])\}$ and
$A_{S}(R[3],S[4])=\{(S[2],S[6]),(S[3],S[5])\}$. In particular,
$(R[2],R[4])$ and $(S[3],S[5])$ are the $R$-parent and $S$-parent
respectively.} \label{F:ancestor}
\end{figure}

An arc is called interior if its start and endpoint are both contained in
either $R$ or $S$ and exterior, otherwise. Let $\prec_1$ be the partial
order $\prec_1$ over the set of interior arcs, given by
\begin{equation}\label{E:partialorder1}
(S[i_1],S[j_1])\prec_1(S[i_2],S[j_2])\quad \Longleftrightarrow\quad
i_2<i_1<j_1<j_2.
\end{equation}
Similarly, let $\prec_2$ denote the partial order over the set of
exterior arcs
\begin{equation}\label{E:porder2}
(S[i_1],R[j_1])\prec_2(S[i_2],R[j_2])\quad \Longleftrightarrow\quad  i_1<i_2,\
j_1<j_2.
\end{equation}
Given an external arc, $(R[i],S[j])$, an interior arc $(R[i_1],R[j_1])$ is
called its $R$-ancestor if $i_1<i<j_1$ and $(S[i_2],S[j_2])$ is the
$S$-ancestor of $(R[i],S[j])$ if $i_2<j<j_2$, respectively.
We call $(R[i],S[j])$ the descendant of $(R[i_1],R[j_1])$ and
$(S[i_2],S[j_2])$ and the sets of $R$-ancestors and $S$-ancestors of
$(R[i],S[j])$ are denoted by $A_{R}(R[i],S[j])$ and $A_{S}(R[i],S[j])$.
The $\prec_1$-minimal $R$-ancestor and $S$-ancestor of $(R[i],S[j])$ are
called its $R$-parent and $S$-parent, see Fig.~\ref{F:ancestor}.
Finally, we call $(R[i_1],R[j_1])$ and $(S[i_2],S[j_2])$ dependent if
they have a common descendant and independent, otherwise.

\begin{figure}[ht]
\centerline{%
\epsfig{file=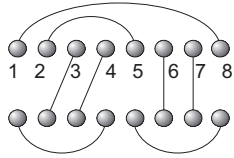,width=0.2\textwidth}\hskip15pt
 }
\caption{\small Subsumed and equivalent arcs: $(R[1],R[8])$ subsumes
$(S[1],S[4])$ and $(S[5],S[8])$. Furthermore, $(R[2],R[5])$ is equivalent
to $(S[1],S[4])$.} \label{F:subsum}
\end{figure}

Suppose $C'=(R',S',I')$ is a subcomplex induced by $R'=R[i_1,j_1]$
and $S'=S[i_2,j_2]$ and suppose furthermore there exists an exterior
arc, $(R[a],S[b])$, with ancestors $(R[i],R[j])$ and
$(S[i'],S[j'])$. The arc $(R[i],R[j])$ is $C'$-subsumed in
$(S[i'],S[j'])$, if for any $(R[k],S[k'])\in I'$ with $i<k<j$, there
exists some $k'$ such that $i'<k'<j'$. In case of $C'=C$, we call
$(R[i],R[j])$ simply ``subsumed'' in $(S[i'],S[j'])$, see
Fig.~\ref{F:subsum}. If $(R[i_1],R[j_1])$ is subsumed in
$(S[i_{2}],S[j_2])$ and vice versa, we call these arcs equivalent.

\begin{figure}[ht]
\centerline{%
\epsfig{file=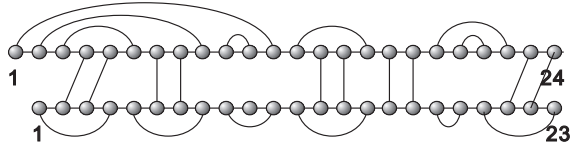,width=0.5\textwidth}\hskip15pt
 }
\caption{\small A joint structure induced by $R[1,24]$ and
$S[1,23]$.} \label{F:joint}
\end{figure}

\begin{figure}[ht]
\centerline{%
\epsfig{file=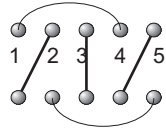,width=0.14\textwidth}\hskip15pt
 }
\caption{\small A zig-zag, generated by $(R[2],S[1])$, $(R[3],S[3])$
and $(R[5],S[4])$.}\label{F:zigzag}
\end{figure}

A joint structure, $J(R[i,j];S[h,\ell],I')=J(i,j;h,\ell)$ is a subcomplex
of $C(R,S,I)$ with the following properties, see Fig.~\ref{F:joint}:\\
$\bullet$ $R$, $S$ are secondary structures\\
$\bullet$ there exist no external pseudoknots, i.e.~if
$(R[i_1],S[j_1]),(R[i_2],S[j_2])\in I'$ where $i_1<i_2$, then
$j_1<j_2$.\\
$\bullet$ there exist no ``zig-zags'', see Fig.\ref{F:zigzag}.
I.e.~if $(R[i_1],R[j_1])$ and $(S[i_2],S[j_2])$ are dependent, then
either $(R[i_1],R[j_1])$ is subsumed by $(S[i_2],S[j_2])$ or vice
versa.
In absence of exterior arcs we refer to a joint structure as a
secondary structure segment, or segment for short. We call
$S[i_1,j_1]$ maximal if there exists no segment, $S[i,j]$,
containing $S[i_1,j_1]$. We remark that the idea of a joint
structure goes back to \cite{Alkan:06} and has also been utilized in
\cite{Backhofen}. One key idea in our approach is to introduce a
specific joint structure, called a tight, which is in some sense a
generalization of the loop. It can be viewed as the transitive
closure of a loop with respect to exterior arcs.

Let $J(a,b;c,d)$ be a fixed joint structure. A joint structure,
$J(i,j;h,\ell)\subset J(a,b;c,d)$ is $J(a,b;c,d)$-tight (or tight in
$J(a,b;c,d)$) if:\\
$\bullet$ there exists at least one exterior arc $(R[i_1],S[j_1])$\\
$\bullet$ for any $(R[i_1],S[j_1])$, we have
\begin{equation}
(A_{R}(R[i_1],S[j_1])\cup A_{S}(R[i_1],S[j_1])) \cap J(a,b;c,d)
 \in J(i,j;h,\ell)
\end{equation}
$\bullet$ $J(i,j;h,\ell)$ is minimal with respect to $\subset$.

\begin{figure}[ht]
\centerline{%
\epsfig{file=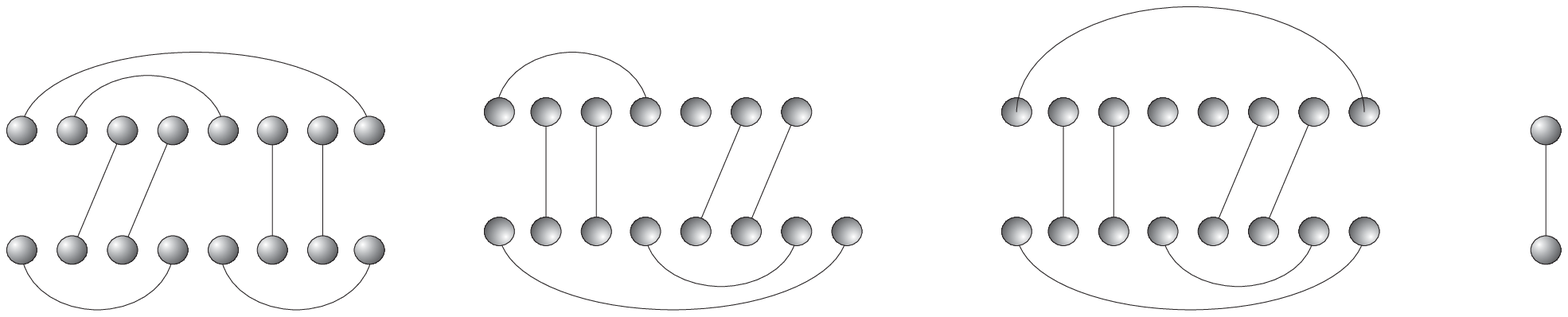,width=0.7\textwidth}\hskip15pt
 }
\caption{\small From left to right: (tjs) of type
$(\bigtriangledown)$, $(\triangle)$, $(\square)$ and $(\circ)$.}
\label{F:typeins}
\end{figure}

\begin{figure}[ht]
\centerline{%
\epsfig{file=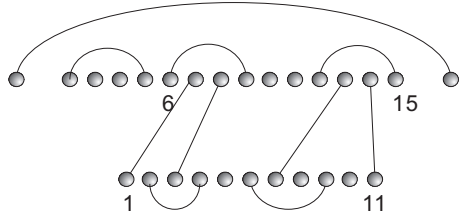,width=0.4\textwidth}\hskip15pt
 }
\caption{\small A double-tight structure $J_{DT}(6,15;1,11)$ in
$J(2,15;1,11)$. Note that the joint structure $J(1,15;1,11)$ itself is a
$\bigtriangledown$-tight.} \label{F:dt}
\end{figure}

Given a tight (tjs), $J(i,j;h,\ell)$, we observe that
neither one of the vertices $i,j,h$ and $\ell$, are start or
endpoint of a segment. In particular, $i,j,h$ and $\ell$ are not
isolated. In combination with the non zig-zag property, we observe
that there are only the following four types of tights
$(\bigtriangledown)$, $(\triangle)$,
$(\square)$ or $(\circ)$, see Fig.\ref{F:typeins}: \\
$(\bigtriangledown)$: $(R[i],R[j])\in J(i,j;h,\ell)$ and
$(S[h],S[\ell])
               \not\in J(i,j;h,\ell)$\\
$(\bigtriangleup)$: $(S[h],S[\ell]) \in J(i,j;h,\ell)$ and
$(R[i],R[j])\not
               \in J(i,j;h,\ell)$\\
$(\square)$: $\{(R[i],R[j]),(S[h],S[\ell])\}\in J(i,j;h,\ell)$\\
$(\circ)$:  $\{(R[i],S[h])\}=J(i,j;h,\ell)$ and $i=j$,
$h=\ell$, i.e.~we have a single interaction.\\
Let $J_{A}(i,j;h,\ell)$ denote a tight structure $J(i,j;h,\ell)$ having
type $\xi$, where $\xi\in A\subset
\{\bigtriangledown,\bigtriangleup,\square,\circ\}$.
In particular, $J_{\xi}(i,j;h,\ell)$ is a tight structure $J(i,j;h,\ell)$
of type $\xi$.
\begin{proposition}\label{P:tight}
Let $J(a,b,c,d)$ be a joint structure, then the following assertions hold:\\
{\rm (a)} if  $J(i,j;h,\ell)$ is tight in $J(a,b,c,d)$, then
$J(i,j;h,\ell)$
          has type $\tau\in \{\bigtriangledown, \triangle, \square,\circ\}$\\
{\rm (b)} any exterior arc is contained in a unique $J(a,b,c,d)$-(tjs)\\
{\rm (c)} $J(a,b,c,d)$ decomposes into a unique sequence of (tjs)
and maximal segments.
\end{proposition}
Suppose we are given two exterior arcs $(R[i_1], S[j_1]),(R[i_2],S[j_2])\in
J(i,j;h,\ell)$. For two $J(i,j;h,\ell)$-tight structures, $J_{T}((R[i_1],
S[j_1]))$, $J_{T}((R[i_2],S[j_2]))$ we set
\begin{equation*}
J_{T}((R[i_1], S[j_1]))=J_{T}((R[i_2],S[j_2]))\quad \Longleftrightarrow\quad
(R[i_1], S[j_1])_1\sim_{J(i,j;h,\ell)} (R[i_2],S[j_2]).
\end{equation*}
Suppose $J_{T}(i,j;r,s)$ is a tight structure where $i\leq a <b\leq
j$ and $r\leq c <d\leq s$. A double-tight structure
$J_{DT}(i,j;r,s)$ in $J_{T}(i,j;r,s)$, is a joint structure
$J(i,j;r,s)$ such that $J(i,j;r,s)\subset J_{T}(i,j;r,s)$ and
\begin{equation}\label{E:TJ2}
J_{DT}(i,j;r,s)=(J_{T}(i,a;r,c),J(a+1,b-1;c+1,d-1),J_{T}(b,j;d,s))
\end{equation}
where $J_{T}(i,a;r,c)$ and $J_{T}(b,j;d,s)$ are
$J(a+1,b-1;c+1,d-1)$-tight structures, see Fig.~\ref{F:dt}.
\begin{figure}[ht]
\centerline{%
\epsfig{file=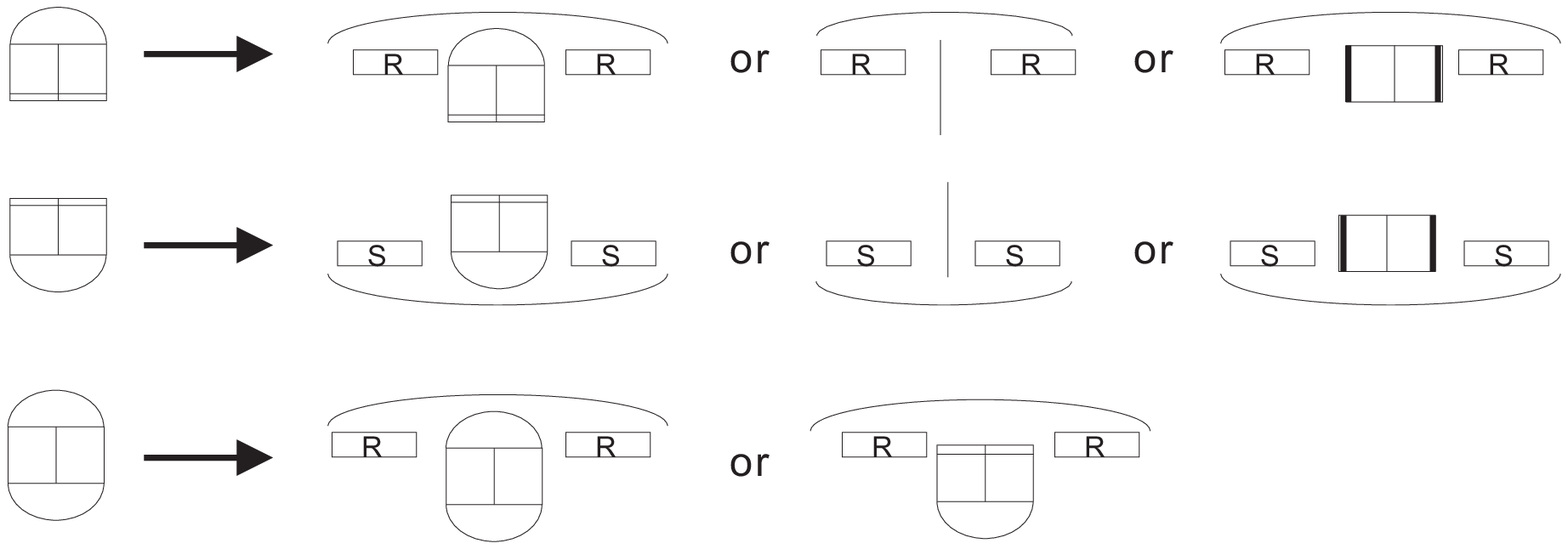,width=0.5\textwidth}\hskip15pt
 }
\caption{\small Decomposing tights: we show how to decompose a
tights of the types $(\bigtriangledown)$, $(\triangle)$ or
$(\square)$ via to Corollary~\ref{C:tu}, Corollary~\ref{C:tu2}
and Corollary~\ref{C:te}.}
\label{F:tiglem}
\end{figure}
\begin{corollary}\label{C:tu}
Let $J_{\bigtriangledown}(i,j;r,s)$ be a tight structure of type
$\bigtriangledown$ and let $\zeta_1=(R[h_1],S[\ell_1])$ and
$\zeta_2=(R[h_2],S[\ell_2])$ be the minimal and maximal exterior
arcs in $J(i,j;r,s)$ and $i+1\leq i_1\leq j_1\leq j-1$. Then
\begin{equation}
\begin{split}
J(i+1,j-1;r,s)=\ \ \ \ \ \ \ \ \ \ \ \ \ \ \ \ \ \ \ \ \ \ \ \ \ \
\ \ \ \ \ \ \ \ \ \ \ \ \ \ \ \ \ \ \ \ \ \ \ \ \ \ \ \ \ \ \ \ \ \ \ \ \ \ \
\ \ \ \ \ \ \ \ \ \ \ \ \ \ \ \ \ \ \ \ \ \ \ \ \ \ \\
\begin{cases}
(R[i+1,i_1-1],J_{\{\bigtriangledown,\circ\}}(i_1,j_1;r,s), R[j_1+1,j-1]), \
& \text{\rm if}
\ \zeta_1\sim_{J(i+1,j-1;r,s)}  \zeta_2 \\
(R[i+1,i_1-1],J_{DT}(i_1,j_1;r,s), R[j_1+1,j-1]) & \text{\rm otherwise,}
\end{cases}
\end{split}
\end{equation}
where $J_{\{\bigtriangledown,\circ\}}(i_1,j_1;r,s)$ denotes a
$J(i+1,j-1;r,s)$-tight of type $\triangle$ or $\circ$.
\end{corollary}
Of course we have
\begin{corollary}\label{C:tu2}
Let $J_{\triangle}(i,j;r,s)$ be a tight structure of type
$\triangle$ and let $\zeta_1=(R[h_1],S[\ell_1])$ and
$\zeta_2=(R[h_2],S[\ell_2])$ be the minimal and maximal exterior
arcs in $J(i,j;r,s)$ and $r+1\leq r_1\leq s_1\leq s-1$. Then
\begin{equation}
\begin{split}
J(i,j;r+1,s-1)=\ \ \ \ \ \ \ \ \ \ \ \ \ \ \ \ \ \ \ \ \ \ \ \ \ \
\ \ \ \ \ \ \ \ \ \ \ \ \ \ \ \ \ \ \ \ \ \ \ \ \ \ \ \ \ \ \ \ \ \ \ \ \ \ \
\ \ \ \ \ \ \ \ \ \ \ \ \ \ \ \ \ \ \ \ \ \ \ \ \ \ \\
\begin{cases}
(S[r+1,r_1-1],J_{\{\triangle,\circ\}}(i,j;r_1,s_1), S[s_1+1,s-1]), \
& \text{\rm if}
\ \zeta_1\sim_{J(i,j;r+1,s-1)}  \zeta_2 \\
(S[r+1,r_1-1],J_{DT}(i,j;r_1,s_1), S[s_1+1,s-1]), & \text{\rm otherwise,}
\end{cases}
\end{split}
\end{equation}
where $J_{\{\triangle,\circ\}}(i_1,j_1;r,s)$ denotes a
$J(i,j;r+1,s-1)$-tight of type $\bigtriangleup$ or $\circ$.
\end{corollary}
\begin{corollary}\label{C:te}
Let $J(i,j;r,s)$ be a tight structure of type $\square$ and set
$i+1\leq i_1\leq j_1\leq j-1$, then $J_{t}(i,j;r,s)$ decomposes as
follows:
\begin{equation}\label{E:TJ3}
J(i+1,j-1;r,s)=(R[i+1,i_1-1],J_{\{\triangle,\square\}}(i_1,j_1;r,s),
R[j_1+1,j-1]),
\end{equation}
where
$J_{\{\triangle,\square\}}(i_1,j_1;r,s)$ denotes a
$J(i+1,j-1;r,s)$-tight of type
$\triangle$ or $\square$.
\end{corollary}
\subsection{Proofs}

\begin{figure}[ht]
\centerline{%
\epsfig{file=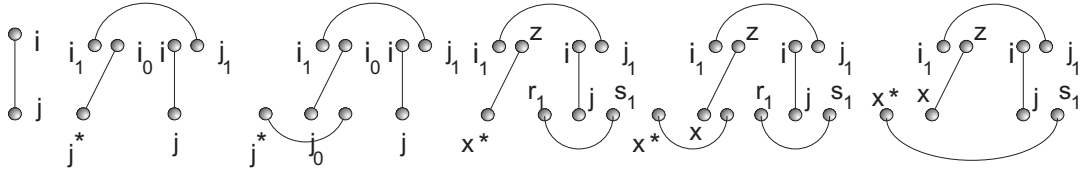,width=0.95\textwidth}\hskip15pt
 }
\caption{\small Illustration of Prop.\ref{P:tight}.} \label{F:prop}
\end{figure}

\textbf{Proof of Proposition \ref{P:tight}}
\begin{proof}
Let $(R[i],S[j])$ be the maximal (rightmost) exterior arc of
$J(a,b,c,d)$. We consider the set of maximal
$(R[i],S[j])$-ancestors, $M$. In case of $M=\varnothing$ we
immediately observe $J(i,j;h,\ell)=(R[i],S[j])$, i.e.~
$J(i,j;h,\ell)$ is of type $\circ$. Suppose next $\vert M\vert =1$.
By symmetry we can, without loss of generality, assume
$M=\{(R[i_1],R[j_1])\}$. Let $(R[i_0],S[j_0])$ the minimal exterior
arc being an descendant of $(R[i_1],R[j_1])$ and let $j_0^*$ denote
either the startpoint of the maximal $(R[i_0],S[j_0])$ $S$-ancestor
or set $j_0^*=j_0$ if no such ancestor exists. Then, by
construction, $J(i_1,j_1;j_0^*,j)$ is tight in $J(a,b,c,d)$.
Finally, in case of $\vert M\vert =2$, i.e.~$M=\{ (R[i_1],R[j_1]),
(S[r_1],S[s_1])\}$. We may, without loss of generality, assume that
$(R[i_1],R[j_1])$ subsumes $(S[r_1],S[s_1])$. Again we consider the
minimal descendant of $(R[i_1],R[j_1])$, $(R[z],S[x])$. Let $x^*$ be
either the startpoint of the maximal $S$-ancestor of $(R[z],S[x])$
or $x^*=x$, otherwise. Then $J(i_1,j_1;x^*,s_1)$ is tight. If
$(R[i_1],R[j_1])$ is equivalent to $S[r_1],S[s_1])$, then
$J(i_1,j_1;r_1,s_1)$ is tight. In the above procedure we have
constructed a (tjs), $J^*$, of type $\tau\in
\{\bigtriangledown,\triangle, \square,\circ\}$ that contains the
maximal exterior $J(a,b,c,d)$-arc. By definition of tight and the
fact that we have noncrossing arcs it follows that any other (tjs)
of $J(a,b,c,d)$ is disjoint to $J^*$. We proceed by considering the
rightmost exterior arc of $J(a,b,c,d)$ that is not contained in
$J^*$, concluding assertion (c) by induction on the number of
exterior arcs of $J(a,b,c,d)$. Since any exterior arc of
$J(a,b,c,d)$ is contained in a unique (tjs) generated by the above
procedure, (b) follows, see Fig.~\ref{F:prop}.
\end{proof}

\textbf{Proof to Corollary \ref{C:tu}}
\begin{proof}
According to Prop. \ref{P:tight}(b), there exist unique
$J(i+1,j-1;r,s)$-tight structures $J(i_1,i_2;r,r_1)$ and
$J(j_2,j_1;s_1,s)$ such that $J(i_1,i_2;r,r_1)=J_{T}(\zeta_1)$ and
$J(j_2,j_1;s_1,s)=J_{T}(\zeta_2)$, respectively.
We have the following two scenarios:
in case of $\zeta_1\sim_{J(i+1,j-1;r,s)}\zeta_2$, i.e.~
$J_{T}(\zeta_1)=J_{T}(\zeta_2)$, we have either $r=s$, in which case
$J(i_1,j_1;r,s)$ is of type $\circ$ and in view of $(S[r],S[s])\not
\in J_{\bigtriangledown}(i,j;r,s)$ $J(i_1,j_1;r,s)$ is of type
$\bigtriangledown$, otherwise.
In case of  $\zeta_1\not \sim_{J(i+1,j-1;r,s)}\zeta_2$, $J(i_1,j_1;r,s)$
is a $J(i+1,j-1;r,s)$-double tight structure.
\end{proof}

\textbf{Proof of Corollary \ref{C:te}}
\begin{proof}
We observe that there exist only one $J(i+1,j-1;r,s)$-tight
structure, since $(S[r],S[s])\in J(i+1,j-1;r,s)$.
We consider the set $M$, consisting of arcs that are equivalent to
$(S[r],S[s])$. According to Prop.~\ref{P:tight}, (c), we have
\begin{equation*}
J(i+1,j-1;r,s)=
\begin{cases}
(R[i+1,i_1-1],J_{\triangle}(i_1,j_1;r,s), R[j_1+1,j-1]) & \text{\rm for}
M=\varnothing\\
(R[i+1,i_1-1],J_{\square}(i_1,j_1;r,s),R[j_1+1,j-1]) & \text{\rm otherwise.}
\end{cases}
\end{equation*}
\end{proof}

\section{Unique decomposition}\label{S:grammar}

\begin{figure}[ht]
\centerline{%
\epsfig{file=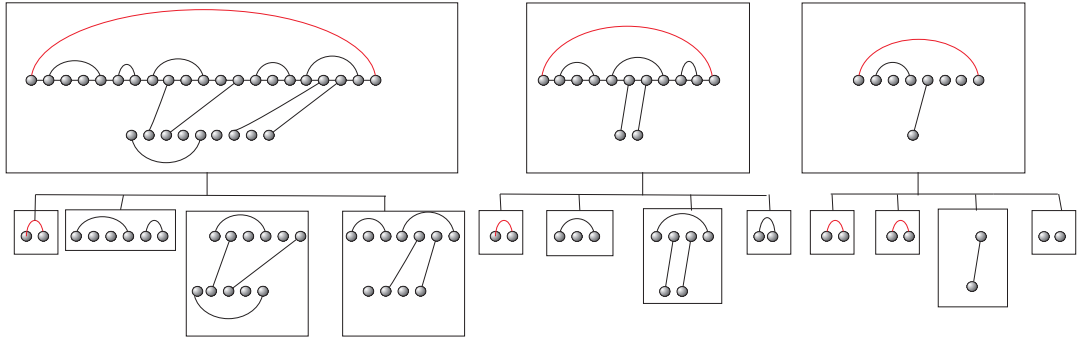,width=0.95\textwidth}\hskip15pt
 }
\caption{\small Illustration of Cor.~\ref{C:tu}.} \label{F:cortu}
\end{figure}

\begin{figure}[ht]
\centerline{%
\epsfig{file=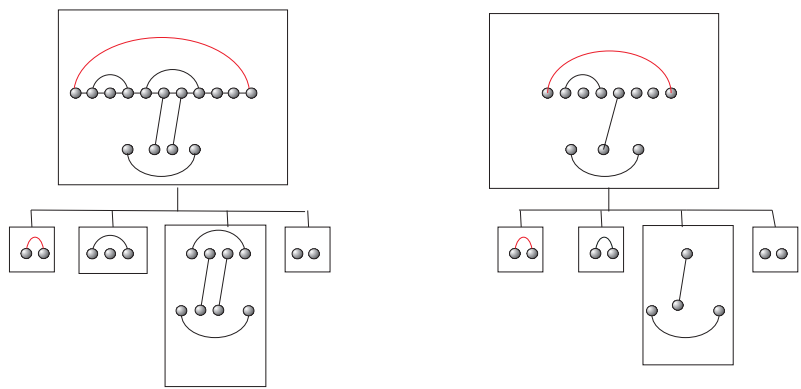,width=0.7\textwidth}\hskip15pt
 }
\caption{\small Illustration of Cor.~\ref{C:te}.} \label{F:corte}
\end{figure}

We showed in Section~\ref{S:comb} via Prop.~\ref{P:tight} that an
arbitrary joint structure uniquely decomposes into a sequence of segments
and tight structures. Via the combinatorial corollaries, Cor.~\ref{C:tu},
Cor.~\ref{C:tu2}
and Cor.~\ref{C:te} we introduced a unique decomposition procedure for tights,
see Fig.~\ref{F:cortu} and Fig.~\ref{F:corte}, below.

\begin{figure}[ht]
\centerline{%
\epsfig{file=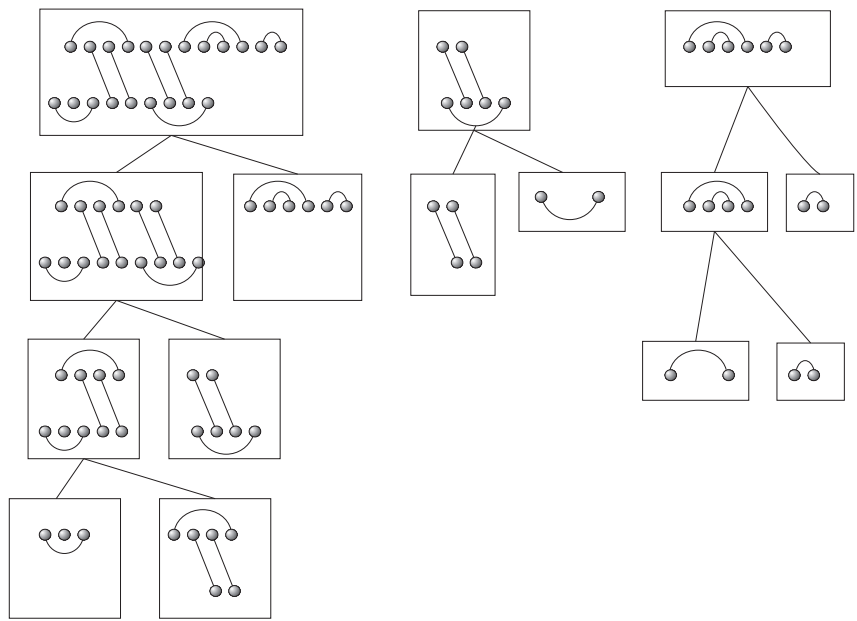,width=0.75\textwidth}\hskip15pt
 }
\caption{\small Illustration of Procedure (a), Procedure (b) and Procedure
(c) for the joint structure $J(1,12,1,8)$.
From left to right we display $T_{a}(1,12;1,8)$,
$T_{b}(5,6;6,9)$ and $T_{c}(R[7,12])$.} \label{F:procedure}
\end{figure}

In this section we give the algorithmic interpretation of the above
results. In the course of our analysis we derive for any joint structure
$J(1,N;1,M)$ a unique decomposition tree via Procedure
(a), (b) and (c), below, see Fig.~\ref{F:procedure}.
Let us begin by giving an interpretation of Prop.~\ref{P:tight}. \\
{\bf Procedure (a):} \\
\underline{input}: a joint structure $\vartheta_{0}=J(i,j;h,\ell)$,
which is not $\vartheta_{0}$-tight or a ms\\
\underline{output}: a unique tree $T_{a}(\vartheta_{0})=(V_{a}(T),E_{a}(T))$\\
Let $i\leq j^* \leq j+1$ and $R[j^*,j]$ be the $\vartheta_{0}$-ms
contain $j$. In particular, $j^*=j+1$ in case of such an ms does not
exist and $j^*=1$ if $R[i,j]$ itself is a ms. Analogously, we define
$S[\ell^*,\ell]$. We construct the
tree $T_{a}(\vartheta_{0})$ recursively as follows: \\
initialization: $V_{a}(T)=\{\vartheta_{0}\}$ and $E_{a}(T)=\varnothing$.\\
(a1): in case of $j^*=j+1$ and $\ell^*=\ell+1$, i.e.~$\vartheta_{0}$
is right-tight, then $\vartheta_{0}$ decomposes  via Prop.~\ref{P:tight}
(b) and (c) into a
$\vartheta_{0}$-tight structure $\vartheta_{1}=J_{\{\bigtriangledown,
\bigtriangleup, \square,\circ\}}(i_1,j;h_1,\ell)$ and a joint structure
$\vartheta_2=J(i,i_1-1;h,h_1-1)$, where $i\leq i_1\leq j$ and
$h\leq h_1\leq \ell$.  Accordingly, we have
\begin{eqnarray}
V_{a}(T)&=& V_{a}(T)\cup \{\vartheta_{1}, \vartheta_{2}\},\\
E_{a}(T)&=& E_{a}(T)\cup
\{(\vartheta_{0},\vartheta_{1}),(\vartheta_{0},\vartheta_{2})\}.
\end{eqnarray}
(a2) otherwise, $\vartheta_{0}$ decomposes into a
$\vartheta_{0}$-right tight structure
$\vartheta_{3}=J_{RT}(i,j^*-1;h,\ell^*-1)$ and two ms
$\vartheta_{4}=R[j^*,j]$, $\vartheta_{5}=S[\ell^*,\ell]$.
Accordingly, we have
\begin{eqnarray}
V_{a}(T)&=&V_{a}(T)\cup \{\vartheta_{3}, \vartheta_{4}, \vartheta_{5}\},\\
E_{a}(T)&=&E_{a}(T)\cup
\{(\vartheta_{0},\vartheta_{3}),(\vartheta_{0},\vartheta_{4}),(\vartheta_{0},
\vartheta_{5})\}.
\end{eqnarray}
We iterate the process until all the leaves of $T_{a}(\vartheta_0)$ are
either $\vartheta_0$-tight structures or $\vartheta_0$-ms.

We proceed by providing an interpretation of Cor.~\ref{C:tu}, Cor.~\ref{C:tu2}
and Cor.~\ref{C:te}.\\
{\bf Procedure (b):} \\
\underline{input}: a tight structure $\vartheta_0=J(i,j;h,\ell)$\\
\underline{output}: a unique tree $T_{b}(\vartheta_{0})=(V_{b}(T),E_{b}(T))$\\
initialization: $V_{b}(T)=\{\vartheta_{0}\}$ and $E_{b}(T)=\varnothing$.\\
We distinguish $J(i,
j;h,\ell)$ by type:\\
$\circ$: do nothing. \\
$\square$: according to Cor.~\ref{C:te},
$\vartheta_0$ decomposes into $\vartheta_{1}=(R[a],R[b])$,
$\vartheta_{2}=R[i+1,i_1-1]$,
$\vartheta_{3}=J_{\square,\bigtriangleup}(i_1,j_1;h,\ell)$ and
$\vartheta_{4}=R[j_1+1,j-1]$, which gives rise to
\begin{eqnarray}
V_{b}(T) &=&V_{a}(T)\cup \{\vartheta_{1}, \vartheta_{2}, \vartheta_{3}
\vartheta_{4}, \vartheta_{5}\},\\
E_{b}(T) &=&E_{a}(T)\cup
\{(\vartheta_{0},\vartheta_{1}),(\vartheta_{0},\vartheta_{2}),(\vartheta_{0},
\vartheta_{3}),(\vartheta_{0},\vartheta_{4}),(\vartheta_{0},\vartheta_{5})\}.
\end{eqnarray}
$\bigtriangledown$: according to Cor \ref{C:tu}, we consider the set of
$J(i+1,j-1;h,\ell)$-tight structures, denoted by $M$. In case of
$\vert M\vert=1$, $J(i+1,j-1;h,\ell)$ decompose into a sequence of a
$J(i+1,j-1;h,\ell)$-tight structure $\vartheta_{6}=J_{\{\bigtriangledown,
\circ\}}(i+1,j-1;h,\ell)$ and two $J(i+1,j-1;h,\ell)$-ms,
$\vartheta_{7}=R[i+1,i_1-1]$ and $\vartheta_{8}=R[j_1+1,j-1]$, where
$i\leq i_1<j_1\leq j$. Accordingly,
\begin{eqnarray}
V_{b}(T)&=&V_{a}(T)\cup \{\vartheta_{1},\vartheta_{6},
\vartheta_{7}, \vartheta_{8}\},\\
E_{b}(T)&=&E_{a}(T)\cup
\{(\vartheta_{0},\vartheta_{1}),(\vartheta_{0},\vartheta_{6}),(\vartheta_{0},
\vartheta_{7}),(\vartheta_{0},\vartheta_{8})\}.
\end{eqnarray}
In case of $\vert M\vert>1$, $J(i+1,j-1;h,\ell)$ decomposes into a sequence
consisting of a $J(i+1,j-1;h,\ell)$-double tight structure
$\vartheta_{9}=J_{DT}(i+1,j-1;h,\ell)$ and two
$J(i+1,j-1;h,\ell)$-ms.~$\vartheta_{7}=R[i+1,i_1-1]$ and
$\vartheta_{8}=R[j_1+1,j-1]$, where $i\leq i_1<j_1\leq j$. Accordingly,
\begin{eqnarray}
V_{b}(T)&=&V_{a}(T)\cup \{\vartheta_{1},\vartheta_{7},
\vartheta_{8}, \vartheta_{9}\},\\
E_{b}(T)&=&E_{a}(T)\cup
\{(\vartheta_{0},\vartheta_{1}),(\vartheta_{0},\vartheta_{7}),(\vartheta_{0},
\vartheta_{8}),(\vartheta_{0},\vartheta_{9})\}.
\end{eqnarray}
Furthermore, let $i_1\leq i_2<j_1$ and $h\leq j_2<\ell$, a
$J(i+1,j-1;h,\ell)$-double tight structure
$\vartheta_{9}=J_{DT}(i+1,j-1;h,\ell)$ decomposes into a
$J(i+1,j-1;h,\ell)$-tight structure
$\vartheta_{10}=J_{\{\bigtriangledown,\circ,\bigtriangleup,\square
\}}(i_1,i_2;h,j_2)$ and a $J(i+1,j-1;h,\ell)$-right tight structure
$\vartheta_{11}=J_{RT}(i_2+1,j_1;j_2+1,\ell)$. I.e.~
\begin{eqnarray}
V_{b}(T)&=&V_{a}(T)\cup \{\vartheta_{10},\vartheta_{11}\},\\
E_{b}(T)&=&E_{a}(T)\cup
\{(\vartheta_{9},\vartheta_{10}),(\vartheta_{9},\vartheta_{11})\}.
\end{eqnarray}
$\bigtriangleup$: analogous to type $\bigtriangledown$ via symmetry.\\
In Fig.~{\ref{F:grammar}} we give an overview of Procedure
(a) and Procedure (b).\\
\begin{figure}[ht]
\centerline{%
\epsfig{file=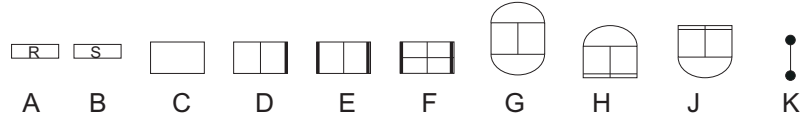,width=0.7\textwidth}\hskip15pt
 }
\caption{\small (A), (B): maximal secondary segments (ms) $R[i,j]$,
$S[r,s]$, (C): joint structure $J(i,j;r,s)$, (D) right-tight
structures $J_{RT}(i,j;r,s)$, (E): double-tight structure
$J_{DT}(i,j;r,s)$, (F): $J_{\bigtriangledown}(i,j;r,s)$, a tight
structure of type $\bigtriangledown$,
    $\bigtriangleup$ or $\square$,
(G): $J_{\square}(i,j;r,s)$, (H): $J_{\bigtriangledown}(i,j;r,s)$,
(J): $J_{\triangle}(i,j;r,s)$ and (K): exterior arc.}
\label{F:graele}
\end{figure}
\begin{figure}[ht]
\centerline{%
\epsfig{file=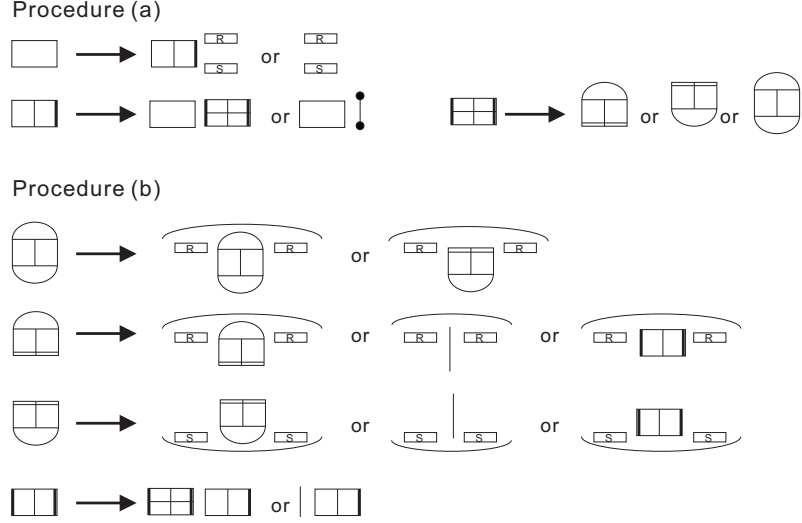,width=0.7\textwidth}\hskip15pt
 }
\caption{\small Illustration of Procedure (a) and Procedure (b),
notations are given via Fig.~{\ref{F:graele} above}.}
\label{F:grammar}
\end{figure}

Finally, we have the wellknown \cite{Waterman:78} secondary structure
loop-decomposition\\
{\bf Procedure $(c)$:}\\
\underline{input}: a secondary structure $\vartheta_0=R[i,j]$\\
\underline{output}: a tree
$T_{c}(\vartheta_{0})=(V_{c}(T),E_{c}(T))$ \\
{initialization}: $V_{b}(T)=\{\vartheta_{0}\}$ and $E_{b}(T)=\varnothing$.\\
We distinguish the following two cases:\\
$(c1)$: in case of $(R[i],R[j])\not\in R[i,j]$, let
$\varnothing_a^b$ denote empty segment in which all the vertices are
isolated. For $1\leq j^*\leq j+1$, let $\varnothing_{j^*}^j$ be the
maximal empty segment that contains $R[j]$. In particular, if $j$ is
not isolated, we have $j^*=j+1$. Let $R^{b}[i_1,j^*-1]$ denote the
segment in which $R[i_1]$ is connected with $R[j^*-1]$. Then
$R[i,j]$ decomposes as follows
$R[i,j]=(\vartheta_1=R[i,i_1-1],\vartheta_2=R^{b}[i_1,j^*-1],
\vartheta_3=\varnothing_{j^*}^j)$ and we set
\begin{eqnarray}
V_{c}(T)&=&V_{c}(T)\cup \{\vartheta_{1},\vartheta_{2},\vartheta_{3}\},\\
E_{c}(T)&=&E_{c}(T)\cup
\{(\vartheta_{0},\vartheta_{1}),(\vartheta_{0},\vartheta_{2}),
(\vartheta_{0},\vartheta_{3})\}.
\end{eqnarray}
\\
$(c2)$: in case of $(R[i],R[j])\in R[i,j]$, i.e.~for $R[i,j]=R^{b}[i,j]$,
we have a decomposition into the pair $(\vartheta_{4}=(R[i],
R[j]),\vartheta_{5}=R[a+1,b-1])$. Accordingly, we have
$V_{c}(T)=V_{c}(T)\cup\{\vartheta_{4},\vartheta_{5}\}$ and
$E_{c}(T)=E_{c}(T)\cup\{(\vartheta_{0},\vartheta_{4}),
(\vartheta_{0},\vartheta_{5})\}$.\\
We iterate (c1) and (c2), until all the leaves in $T$ are
either isolated segments or single arcs.

\begin{figure}[ht]
\centerline{%
\epsfig{file=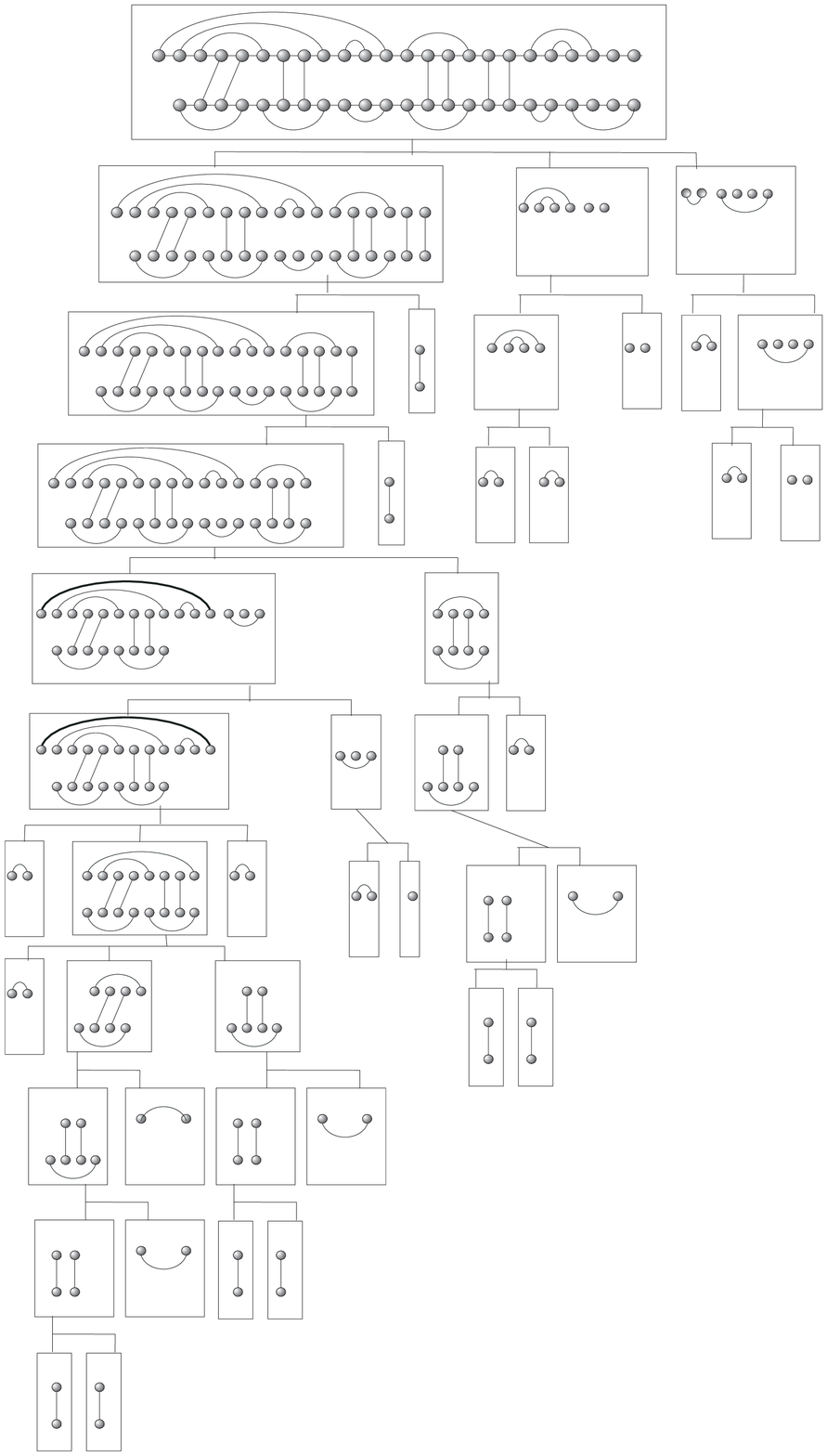,width=0.75 \textwidth}\hskip15pt
 }
\caption{\small The decomposition tree $T(1,N;1,M)$.}
\label{F:decomtree}
\end{figure}

For any joint structure, $J(1,N;1,M)$, we can now construct a tree, with
root $J(1,N;1,M)$ and whose vertices are specific
subgraphs of $J(1,N;1,M)$. The latter are obtained by successive application
of Procedure (a), (b) and (c), see Fig.~\ref{F:decomtree}.
To be precise, let $H$ be the graph rooted in $J(1,N;1,M)$ defined
inductively as follows: for the induction basis for fixed
$J(1,N;1,M)$ only one, Procedure (a), (b) or (c) applies. Procedure
(a), (b) or (c) generates the (procedure-specific, nontrivial)
subtrees, $T_{a}$, $T_{b}$ and $T_{c}$. Suppose
$\vartheta_{\dagger}$ is a leaf of $T$ that has been constructed via
Procedure (a), (b) or (c). As in case of the induction basis, each
such leaf is input for exactly one procedure, which in turn
generates a corresponding subtree. Prop.~\ref{P:tight},
Cor.~\ref{C:tu}, Cor.~\ref{C:tu2} and Cor.~\ref{C:te} imply that
$H$ itself is a tree.
We denote this decomposition tree by $T(1,N;1,M)$, see
Fig.~\ref{F:decomtree}. Accordingly, we have proved

{\bf Observation 1.} {\it For any joint structure, $J(1,N;1,M)$, there
exists a unique decomposition tree, $T(1,N;1,M)$, whose leafs are
either interior or exterior $J(1,N;1,M)$-arcs or isolated segments. }

As we shall see in Section~\ref{S:prob}, the decomposition tree
plays a key role for the calculation of the base pairing
probabilities. To be precise, given a joint structure,
$J(i,j;h,\ell)$, let $T_J(1,N;1,M)$ be the decomposition tree of
$J(1,N;1,M)$ and let $\Sigma_{0}=\{J(1,N;1,M) \mid J(i,j;h,\ell)\in
T_J(1,N;1,M)\}$. Then the probability of $J(i,j;h,\ell)$, denoted by
$\mathbb{P}(i,j;h,\ell)$, is given by
\begin{equation}\label{E:Jprob}
\mathbb{P}(i,j;h,\ell)=\sum_{J(1,N;1,M)\in
\Sigma_{0}}\mathbb{P}(J(1,N;1,M))
\end{equation}
and furthermore

\begin{figure}[ht]
\centerline{%
\epsfig{file=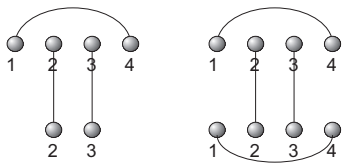,width=0.3\textwidth}\hskip12pt
 }
\caption{\small $J(1,4;2,3)$ has the property that
$J(1,4;2,3)\subset J(1,4;1,4)$ but $J(1,4;2,3)\not\in
T_J(1,4;1,4)$.} \label{F:different}
\end{figure}

{\bf Observation 2.} {\it In general $J(i,j;h,\ell)\subset
J(1,N;1,M)$ is not equivalent to $J(i,j;h,\ell)\in T_J(1,N;1,M)$, see
Fig.~\ref{F:different}. However, in case of secondary structures,
i.e.~$J(i,j;h,\ell)= (R[i],R[j])$, we have
\begin{equation}
(R[i],R[j])\subset J(1,N;1,M)\Leftrightarrow (R[i],R[j])\in
T_J(1,N;1,M).
\end{equation}}

\section{From the decomposition tree to the partition function}
\label{S:partition}

We discussed in the introduction the concept of the loop-based
partition function of RNA secondary structures due to McCaskill
\cite{McCaskill}. We observed there that the key property for its
derivation is the unique decomposition into substructures and their
recursive analysis. For instance, suppose we are given a tight of
type $\bigtriangledown$ from which we remove, by virtue of
Cor.~\ref{C:tu}, its outer arc. For this purpose, the context of the
latter, i.e.~its particular arc-configuration has to be taken into
account. However, once the unique decomposition is established, the
existence of specific subclasses of joint structures allowing for
the dynamic programming of the partition function follows. We remark
that the particular choice of the latter may not be unique.

\begin{figure}[ht]
\centerline{%
\epsfig{file=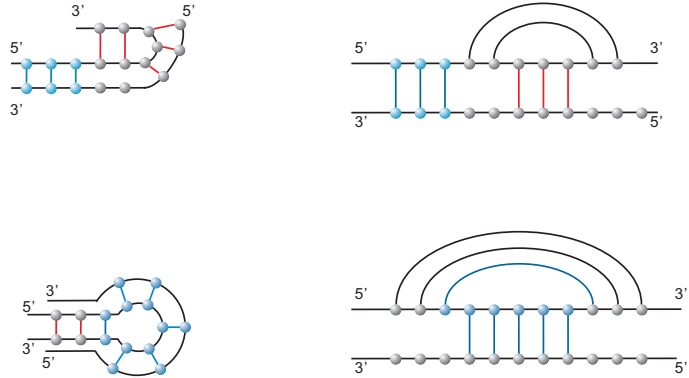,width=0.6 \textwidth}\hskip15pt
 }
\caption{\small The two new loop types: the hybrid (top) and the
kissing loop (bottom).}
\label{F:kissingloop}
\end{figure}

The first step is to extend the standard loop-energy model for
secondary structures by introducing two new loop-types due to Chitsaz
{\it et al.} \cite{Backhofen}: the kissing loop and the hybrid,
see Figure~\ref{F:kissingloop}.

\subsection{Loops}
Having discussed the standard loop types of secondary structures in
Section~\ref{S:Introduction}, we proceed now by introducing the loops that
contain exterior arcs.\\
{\bf (4)} a {\it hybrid}-loop (${\sf Hy}$)is a sequence
$((R[i_1],S[j_1]),\dots,([R[i_s],S[j_s]))$, where $s\ge 2$ and
$(i_r,j_r)$ is nested in $(i_1,j_1)$ such that $R[i_h+1,i_{h+1}-1]=
\varnothing_{i_h+1}^{i_{h+1}-1}$ and $S[j_h+1,j_{h+1}-1]=
\varnothing_{j_h+1}^{j_{h+1}-1}$.\\
{\bf (5)} a {\it kissing}-loop (${\sf K}$) is either a pair,
$\left((R[i],R[j]),R[i+1,j-1]\right)$, such that there exists at
least one $(R[i],R[j])$-child, $(R[i_1],S[j_1])$ where $i<i_1<j$ or
a pair $\left((S[i],S[j]),S[i+1,j-1]\right)$, with
$(R[i],R[j])$-child $(R[i_1],S[j_1])$ and $i<j_1<j$.\\

The arguments of Prop.~\ref{P:tight}, Cor.~\ref{C:tu},
Cor.~\ref{C:tu2} and Cor.~\ref{C:te} imply that each joint structure
can uniquely be decomposed into a sequence of loops--a necessary and
sufficient condition for the mfe-folding of joint structures. As we
shall see in the next section, the unique decomposition {\it and}
the particular choice of loops give rise to specific subclasses via
which the partition function can be recursively expressed.
Furthermore, following \cite{Bernhart}, we allow for an initiation
energy, i.e.~each hybrid loop is given an energy penalty of
$\sigma_{0}$. In addition, we allow for a scaling, $0<\sigma\leq 1$,
of the energy contribution of each hybrid loop. As default we set
$\sigma_{0}=0$, $\sigma=1$.

\subsection{Case studies}
Consider a joint structure $J(i,j;h,\ell)\in T(J(1,N;1,M))$. For the
purpose of assigning an energy to a substructure, we have to
distinguish substructures by their ``outer" loop type, see Case $1$
as well as Fig.~\ref{F:stand} and Fig~\ref{F:kissingloop}. To convey the
key ideas we shall restrict our analysis to three case studies.
\begin{figure}[ht]
\centerline{%
\epsfig{file=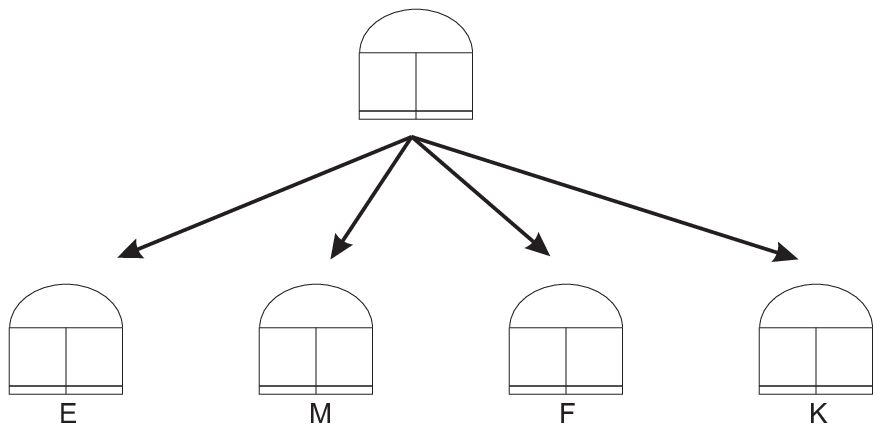,width=0.4\textwidth}\hskip12pt
 }
\caption{\small Context dependency: the labels ``${\sf E},{\sf
M},{\sf F},{\sf K}$'' are defined in Case $1$. We display from left
to right $J^{{\sf E}}_{\bigtriangledown}(i,j;h,\ell)$, $J^{\sf
{M}}_{\bigtriangledown}(i,j;h,\ell)$,
$J^{\sf{F}}_{\bigtriangledown}(i,j;h,\ell)$ and
$J^{\sf{K}}_{\bigtriangledown}(i,j;h,\ell)$,
respectively.}\label{F:tdtype}
\end{figure}

Given a joint structure $J(i,j;h,\ell)\in T(J(1,N;1,M))$, we set
$M_{R}(i,j)=\{(R[i_1],R[j_1])\mid i_1<i\leq j<j_1\}$ and
$M_{S}(h,\ell)=\{(S[i_1],S[j_1])\mid i_1<h\leq \ell<j_1\}$.

{\it Case $1$.}
Suppose we are given a tight structure $J_{\bigtriangledown}(i,j;h,\ell)$.
In case of $M_{S}(h,\ell)=\varnothing$, we call $S[h,\ell]$ external and use
the notation $J^{{\sf E}}_{\bigtriangledown}(i,j;h,\ell)$.
Otherwise, let $(S[i_0],S[j_{0}])$ be the minimal element of $M_{S}(h,\ell)$.
We denote the type of the loop including $(S[i_0],S[j_{0}])$, by $\xi$. In
case of $\xi={\sf M}$, we use the notation
$J^{{\sf M}}_{\bigtriangledown}(i,j;h,\ell)$. Otherwise, in case of
$\xi={\sf K}$, we write
$J^{{\sf K}}_{\bigtriangledown}(i,j;h,\ell)$ or
$J^{{\sf F}}_{\bigtriangledown}(i,j;h,\ell)$ depending on whether or not
$J_{\bigtriangledown}(i,j;h,\ell)$ contains the child of
$(S[i_0],S[j_{0}])$, see Fig.~\ref{F:tdtype}.

\begin{figure}[ht]
\centerline{%
\epsfig{file=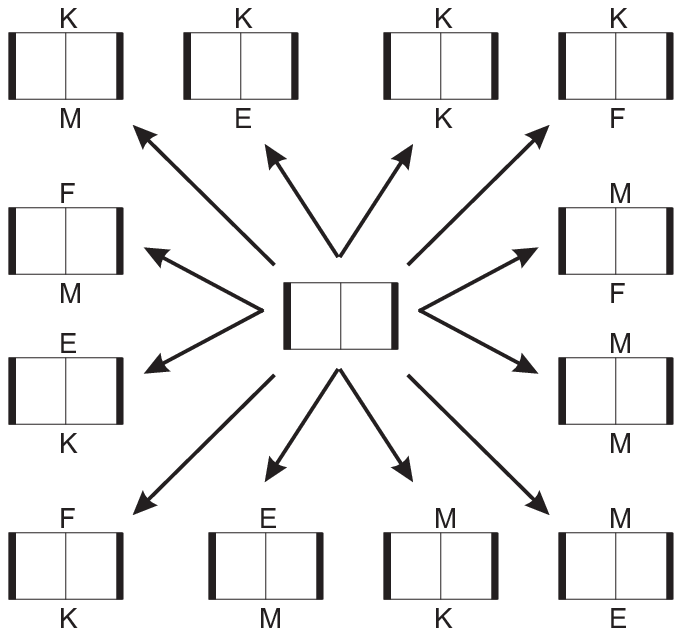,width=0.4\textwidth}\hskip12pt
 }
\caption{\small The twelve subclasses of $J_{DT}(i,j;h,\ell)$ as discussed in
Case $2$.}
\label{F:dttype}
\end{figure}

{\it Case $2$.} Suppose we are given a double-tight structure,
$J_{DT}(i,j;h,\ell)$. Then we arrive at the twelve subclasses
presented in Figure~\ref{F:dttype}. Indeed, according to
Cor.~\ref{C:tu}, there does not exist any $J_{DT}^{{\sf E},{\sf
E}}(i,j;h,\ell)$, i.e.~$M_{R}(i,j)\cup M_{S}(h,\ell)\neq
\varnothing$. Without loss of generality, we may assume that
$M_{R}(i,j)\neq \varnothing$ and that $(R[i_{1}],R[j_1])\in
M_{R}(i,j)$ is minimal. In case of $M_{S}(h,\ell)=\varnothing$, we
use the notation $J_{DT}^{Y,{\sf E}}(i,j;h,\ell)$, where $Y$ is the
loop type formed by $(R[i_{1}],R[j_1])$ and $R[i_1+1,j_1-1]$.
Otherwise, we have $M_{S}(h,\ell)\neq\varnothing$. Let
$(S[i_{2}],S[j_2])$ be the minimal element. In this case we use the
notation $J_{DT}^{Y_1,Y_2}(i,j;h,\ell)$, where $Y_1$ and $Y_2$ are
the loop-types formed by $(R[i_{1}],R[j_1])$, $R[i_1+1,j_1-1]$ and
$(S[i_{2}],S[j_2])$, $S[i_2+1,j_2-1]$, respectively, see
Fig.~\ref{F:dttype}.

\begin{figure}[ht]
\centerline{%
\epsfig{file=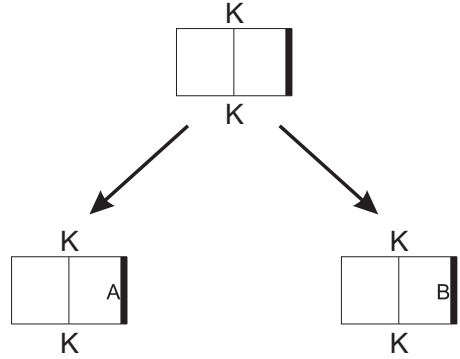,width=0.4\textwidth}\hskip12pt
 }
\caption{\small The four subclasses of $J_{RT}^{{\sf K},{\sf
K}}(i,j;h,\ell)$, see Case $3$ for details.} \label{F:dtuu}
\end{figure}

{\it Case $3$.} In case of a right-tight structure, $J_{RT}^{{\sf
K},\sf{K}}(i,j;h,\ell)$, we obtain four subclasses. In case of
$(R[j],S[\ell])\in J_{RT}^{{\sf K},{\sf K}}(i,j;h,\ell)$, we say
$J_{RT}^{{\sf K},{\sf K}}(i,j;h,\ell)$ is $(rB)$ and $(rA)$,
otherwise. Let $(R[i_{1}],S[j_1])$ denote the minimal exterior arc
in $J_{RT}^{{\sf K},{\sf K}}(i,j;h,\ell)$. According to
Prop.~\ref{P:tight}, there exists a unique $J(i,j;h,\ell)$-tight
structure $J_{T}(R[i_{1}],S[j_1])$, such that $(R[i_{1}],S[j_1])\in
J_{T}(R[i_{1}],S[j_1])$. In case of $J_{T}(R[i_{1}],S[j_1])$ is of
type $\circ$, i.e. $(R[i_{1}],S[j_1])$ itself and
$R[i,i_1]=\varnothing_{i}^{i_1}$,
$S[h,j_{1}]=\varnothing_{h}^{j_1}$, we say $J_{RT}^{{\sf K},{\sf
K}}(i,j;h,\ell)$ is $(lB)$ and $(lA)$, otherwise. We use the
notation $J_{RT}^{{\sf K},{\sf K},Y_1,Y_2}(i,j;h,\ell)$, if
$J_{RT}^{{\sf K},{\sf K}}(i,j;h,\ell)$ is $(lY_1)$ and $(rY_2)$,
respectively, see Fig.~\ref{F:dtuu}.

\subsection{The partition function}
\begin{figure}[ht]
\centerline{%
\epsfig{file=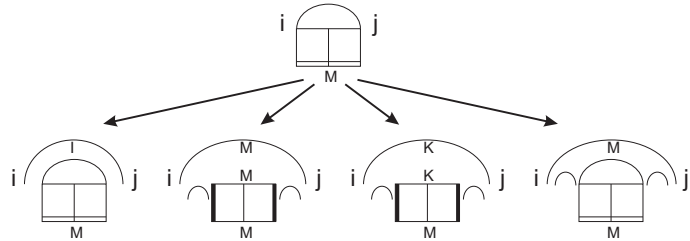,width=0.6\textwidth}\hskip12pt
 }
\caption{\small The four decompositions of
$J^{\sf{M}}_{\bigtriangledown}(i,j;r,s)$ via Procedure (b), denoted
by $I_{1}$, $I_2$, $I_{3}$ and $I_{4}$, from left to right,
respectively.} \label{F:parexp1}
\end{figure}

In the previous section we discussed specific subclasses of joint structures.
They were designed to facilitate the recursive construction of the
partition function. The purpose of this section is to showcase the
respective recursions induced by these classes.\\
{\it Case $1$:} $J^{\sf{M}}_{\bigtriangledown}(i,j;r,s)$. According
to Cor.~\ref{C:tu}, we have three cases:
$J_{\bigtriangledown}(i,j;r,s)$ decomposes into either a
$J(i-1,j+1;r,s)$-tight structure of type $\kappa$, where $\kappa\in
\{\bigtriangledown,\circ\}$ or a $J(i-1,j+1;r,s)$-double tight
structure and a ms. By definition of
$J^{\sf{M}}_{\bigtriangledown}(i,j;r,s)$, the case of a
$J(i-1,j+1;r,s)$-tight structure of type $\circ$ is impossible.
Considering the type of the loop including $(R[i],R[j])$ and
$R[i+1,j-1]$, we arrive exactly at the four cases, denoted by
$I_{1}$, $I_2$, $I_{3}$ and $I_{4}$, from left to right, displayed
in Fig.~\ref{F:parexp1}.

\begin{figure}[ht]
\centerline{%
\epsfig{file=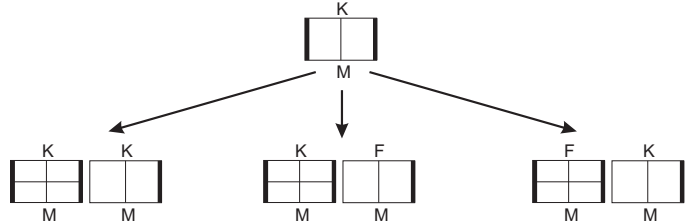,width=0.6\textwidth}\hskip12pt
 }
\caption{\small The decomposition of $J_{DT}^{K,M}(i,j;h,\ell)$ via
Procedure (b). The corresponding three cases are labeled by from
left to right by $J_{1}$, $J_2$ and $J_{3}$, respectively.}
\label{F:parexp2}
\end{figure}

\begin{figure}[ht]
\centerline{%
\epsfig{file=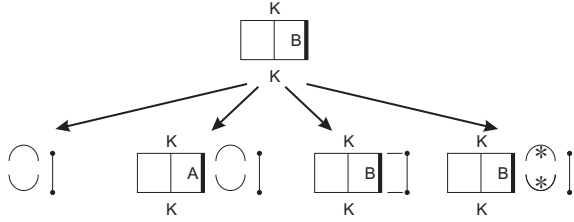,width=0.5\textwidth}\hskip12pt
 }
\caption{\small The four decomposition scenarios of
$J^{{\sf{K},\sf{K}},B,B}_{DT}(i,j;r,s)$ via Procedure (b). We denote
the corresponding cases from left to right by $L_{1}$, $L_2$ $L_3$
and $L_{4}$, respectively.} \label{F:parexp3}
\end{figure}

Let $i<h<\ell<j$. According to the recurrences displayed in
Fig.~\ref{F:parexp1}, the partition function satisfies for
$J^{\sf{M}}_{\bigtriangledown}(i,j;r,s)$ the following recursion:\\
\begin{equation}
Q^{\sf{M}}_{\bigtriangledown}(i,j;r,s)=\sum_{h,\ell}
(Q(I_{1})+Q(I_{2})+Q(I_{3})+Q(I_{4})),
\end{equation}
where
\begin{eqnarray*}
Q(I_{1})&=&Q_{\bigtriangledown}^{\sf{M}}
(h,\ell;r,s)e^{-{\sf Int}(i,j;h,\ell)/kT},\\
Q(I_{2})&=&Q_{DT}^{\sf{M},\sf{M}}(h,\ell;r,s)e^{-(\alpha_1+\alpha_2)/kT}
(e^{-(j-\ell-1)\alpha_3/kT}
+Q^{\sf{m}}(i+1,h-1))\\
&&\times(e^{-(h-i-1)\alpha_3/kT}+Q^{\sf{m}}(\ell+1,j-1)),\\
Q(I_{3})&=&Q_{DT}^{\sf{K},\sf{M}}(h,\ell;r,s)e^{-(\beta_1+\beta_2)/kT}
(e^{-(j-\ell-1)\beta_3/kT}
+Q^{\sf{k}}(i+1,h-1))\\
&&\times(e^{-(h-i-1)\beta_3/kT}+Q^{{\sf k}}(\ell+1,j-1)),\\
Q(I_{4})&=&Q_{\bigtriangledown}^{\sf{M}}(h,\ell;r,s)e^{-(\alpha_1+2\alpha_2)/kT}
(Q^{-(j-\ell-1)\alpha_3/kT}Q^{\sf{m}}(i+1,h-1)\\
& &+e^{-(h-i-1)\alpha_3/kT}
Q^{\sf{m}}(\ell+1,j-1))+Q^{\sf{m}}(\ell+1,j-1)Q^{\sf{m}}(i+1,h-1)).\\
\end{eqnarray*}
{\it Case 2:} $J_{DT}^{{\sf K},{\sf M}}(i,j;h,\ell)$. According to
Procedure (b), a double tight structure decomposes into a
$J(i,j;h,\ell)$-tight structure, $J(i,i_1;h,h_1)$ and a
$J(i,j;h,\ell)$-right tight structure, $J(i_1+1,j;h_1+1,\ell)$. We
observe that the type of the outer loop of $S[h,h_1]$ and
$S[h_1+1,\ell]$ coincides with that of $S[h,\ell]$, i.e.~${\sf M}$.
Analogously, the outer loop of $R[i,i_1]$ and $R[i_1+1,j]$, denoted
by $(R[i_{0},j_{0}])$, is of type ${\sf K}$. Furthermore, at least
one of the substructures $R[i,i_1]$ and $R[i_1+1,j]$ contain the
child of $(R[i_{0},j_{0}])$. Consequently we arrive at the three
scenarios labeled by from left to right by $J_{1}$, $J_2$ and
$J_{3}$ displayed in Fig.~\ref{F:parexp2}.
Setting
\begin{equation*}
Q_{\bigtriangleup,\bigtriangledown, \square}^{{\sf K},{\sf
M}}(i,i_1;r,j_1)=Q_{\bigtriangleup}^{{\sf K},{\sf M}}(i,i_1;r,j_1)
+Q_{\bigtriangledown}^{{\sf K},{\sf
M}}(i,i_1;r,j_1)+Q_{\square}^{{\sf K},{\sf M}}(i,i_1;r,j_1),
\end{equation*}
the recursion of the partition function for $J_{DT}^{{\sf K},{\sf
M}}(i,j;h,\ell)$
is given by:\\
\begin{equation}
Q_{DT}^{{\sf K},{\sf M}}(i,j;r,s)=\sum_{i_1,j_1}(Q(J_{1})+Q(J_{2})
+Q(J_{3})),
\end{equation}
where
\begin{eqnarray*}
Q(J_{1})&=&Q_{\bigtriangleup,\bigtriangledown,
\square}^{{\sf K},{\sf M}}(i,i_1;r,j_1)
Q_{DT}^{{\sf K},{\sf M}}(i_1+1;j_1+1,s)\\
Q(J_{2})&=&Q_{\bigtriangleup,\bigtriangledown,
\square}^{{\sf F},{\sf M}}(i,i_1;r,j_1)
Q_{DT}^{{\sf K},{\sf M}}(i_1+1;j_1+1,s))\\
Q(J_{3})&=&Q_{\bigtriangleup,\bigtriangledown, \square}^{{\sf
K},{\sf M}}(i,i_1;r,j_1)Q_{DT}^{{\sf F},{\sf M}}(i_1+1;j_1+1,s).
\end{eqnarray*}
{\it Case 3:} $J_{DT}^{{\sf{K},\sf{K}},B,B}(i,j;h,\ell)$. By
definition of $J_{DT}^{{\sf{K},\sf{K}},B,B}(i,j;h,\ell)$, we have
$(R[j],S[\ell])\in J(i,j;h,\ell)$. We consider the set of exterior
arcs in $J(i,j-1;r,\ell-1)$, $W$. In case of $W=\varnothing$,
$J_{DT}^{{\sf{K},\sf{K}},B,B}(i,j;h,\ell)$ decomposes into
$R[i,j-1]$, $S[h,\ell-1]$ and $(R[j],S[\ell])$. This is the leftmost
(first) case ($L_1$) displayed in Fig.~\ref{F:parexp3}. Otherwise, let
$(R[i_{1}],S[j_1])$ denote the maximal exterior arc in
$J(i,j-1;r,\ell-1)$. We consider the unique $J(i,j;r,\ell)$-tight
structure which contains $(R[i_{1}],S[j_1])$, denoted by
$J_{T}(R[i_{1}],S[j_1])$. If $J_{T}(R[i_{1}],S[j_1])$ has not type
$\circ$, we have the second case ($L_2$) displayed in Fig.~\ref{F:parexp3}.
Otherwise, depending on whether or not
$R[j_1+1,j-1]=\varnothing_{j_1+1}^{j-1}$ and
$S[h_1+1,\ell-1]=\varnothing_{h_1+1}^{\ell-1}$, we have the third ($L_3$)
and fourth case ($L_4$), displayed in Fig.~\ref{F:parexp3}.
Consequently, we arrive at:\\
\begin{equation}
Q_{DT}^{{\sf{K},\sf{K}},B,B}(i,j;h,\ell)=\sum_{j_1,h_1}(Q(L_{1})+Q(L_{2})
+Q(L_{3})+Q(L_{4})),
\end{equation}
where $Q(L_1)=Q_{s}^{\sf k}(i,j-1)Q_{s}^{\sf k}(h,\ell-1)$ and
\begin{eqnarray*}
Q(L_{2})&=&Q_{RT}^{{\sf{K},\sf{K}},B,A}(i,j_1;h,h_1)
(Q^{\sf{k}}(j_1+1,j-1)+Q^{\sf{k}}(h_1+1,\ell-1)
+Q^{\sf{k}}(j_1+1,j-1)Q^{\sf{k}}(h_1+1,\ell-1))\\
Q(L_{3})&=&Q_{RT}^{{\sf{K},\sf{K}},B,B}(i,j_1;h,h_1)e^{-\sigma_{0}+\sigma
{\sf{Int}}(j,\ell;j_1,h_1)}\\
Q(L_{4})&=&Q_{RT}^{{\sf{K},\sf{K}},B,B}(i,j_1;h,h_1)
(Q^{\sf{k}}(j_1+1,j-1)+Q^{\sf{k}}(h_1+1,\ell-1)
+Q^{\sf{k}}(j_1+1,j-1)Q^{\sf{k}}(h_1+1,\ell-1))
\end{eqnarray*}

\section{Base pairing probabilities}\label{S:prob}

We have seen in Section~\ref{S:grammar} that the probability of a
joint structure, $J(1,N;1,M)$, is given by
\begin{equation}
\mathbb{P}(J(1,N;1,M))=\frac{1}{Q^{I}}e^{-F(J(1,N;1,M))/kT},
\end{equation}
where $Q^I=\sum_{J(1,N;1,M)} e^{-F(J(1,N;1,M))/kT}$. In this
section, we shall calculate the base pair probabilities (BPP) for
interior and exterior arcs. The key idea is here to associate the
probability of specific substructures contained in the decomposition
tree. In other words, a term $Q^{Y_1,Y_2,Y_3,Y_4}_{\xi}(i,j;h,\ell)$
in the recursive calculation of the partition function gives rise
to the probability $\mathbb{P}^{Y_1,Y_2,Y_3,Y_4}_{\xi}(i,j;h,\ell)$.
For instance, $\mathbb{P}^{{\sf M},{\sf K},A,B}_{RT}(i,j;h,\ell)$
is, by construction, the sum over all the probabilities of joint
structures $J(1,N;1,M)$ such that $J(i,j;h,\ell)$ is contained in
$T(J(1,N;1,M))$ and $J(i,j;h,\ell)= {J}^{{\sf M},{\sf
K},A,B}_{RT}(i,j;h,\ell)$. We remark that the above observations
reduce the computation of the BPP to a trace-back routine in the
decomposition tree, constructed in Section~\ref{S:grammar}.

The basic strategy can be sketched as follows: \\
{\bf (a)} derive from the recursion of the partition function the
corresponding recursion of the probabilities\\
{\bf (b)} partition the substructures according to their respective
contribution to the partition function\\
{\bf (c)} for each subclass, recursively calculate the probability
of substructures via tracing back the decomposition tree.

We recall that $\Sigma_{0}=\{J(1,N;1,M)\mid J(i,j;h,\ell) \in
T(J(1,N;1,M))\}$. The probability  $\mathbb{P}(i,j;h,\ell)$ is given
by
\begin{equation}
\mathbb{P}(i,j;h,\ell)=\sum_{J(1,N;1,M)\in
\Sigma_{0}}\mathbb{P}(J(1,N;1,M)).
\end{equation}
We accordingly set
\begin{equation}\label{E:Jprob1}
\mathbb{P}^{Y_1,Y_2,Y_3,Y_4}_{\xi}(i,j;h,\ell)=\sum_{J(1,N;1,M)\in
\Lambda_{0}}\mathbb{P}(J(1,N;1,M)),
\end{equation}
where $\Lambda_{0}=\{J(1,N;1,M)\mid J(i,j;h,\ell)\in T(J(1,N;1,M)),
J(i,j;h,\ell)\in J^{Y_1,Y_2,Y_3,Y_4}_{\xi}(i,j;h,\ell)\}$.

\subsection{Base pairing probabilities for RNA secondary structures}
\begin{figure}[ht]
\centerline{%
\epsfig{file=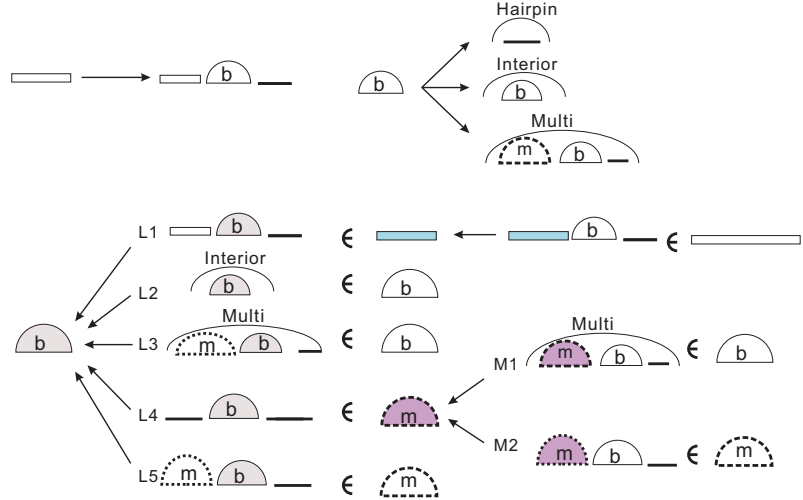,width=0.7 \textwidth}\hskip15pt
 }
\caption{\small Tracing back: for a parent of $R^{b}(i,j)$ we have
according to Procedure (c), five cases, labeled from top to bottom
by $L_1$, $L_2$, $L_3$, $L_4$ and $L_5$. For a parent of
$R^{m}(i,j)$ there are two cases, denoted by $M_1$ and
$M_2$.}\label{F:secinverse}
\end{figure}

In order to illustrate the concept, let us put the calculation of
the BPP for secondary structures into the context of our
backtracking routine. Given a secondary structure $R$ of length $N$,
the probability of $R$ is given by
$\mathbb{P}(R)=\frac{1}{Q}e^{-F(R)/kT}$. In order to calculate the
probability of $R[i]$ being connected to $R[j]$ in the equilibrium
ensemble of structures, $\mathbb{P}(i_R,j_R)$, the first objective
is to express the probability of this base pair into a sum of
probabilities of substructures. Let $T(R[1,N])$ be the decomposition
tree of a particular secondary structure $R[1,N]$ via Procedure (c)
and $\Omega(i_{R},j_{R})=\{S\mid (R[i],R[j])\in T(S)\}$. We remark
that $\Omega(i_{R},j_{R})$ coincides set of secondary structure such
that $R[i]$ is bound with $R[j]$, see Section~\ref{S:grammar},
Observation 2. Then we have
\begin{equation}
\mathbb{P}(i_{R},j_{R})=\frac{\sum_{S\in\Omega(i_R,j_R)}Q(S)}{Q}.
\end{equation}
Let $R^{b}(i,j)$ denote the set of segments $R[i,j]$ in which $R[i]$
is connected with $R[j]$ and $R[i,j]\in T(R[1,N])$. By construction,
$\mathbb{P}^{b}(i_{R},j_{R})$ is the probability of $R^{b}(i,j)$.
According to Procedure (c), we have
$\mathbb{P}(i_{R},j_{R})=\mathbb{P}^{b}(i_{R},j_{R})$ since
$(R[i],R[j])\in T(J(1,N;1,M))$ if and only if the parent of
$(R[i],R[j])$ in the decomposition tree belongs to $R^{b}(i,j)$.
Therefore the problem is reduced to the calculation of
$\mathbb{P}^{b}(i_{R},j_{R})$. Inspection of Procedure (c) shows,
that for the parent of an element of $R^{b}(i,j)$ we have to
distinguish the five cases displayed in Fig.~\ref{F:secinverse}.
Let $R^{m}(i,j)$ denote the set of segments $R[i,j]\in T(R[1,N])$
such that $R[i,j]\neq \varnothing_{i}^{j}$, where the outer loop has
type ${\sf M}$. Let $R^{s}(i,j)$ denote the set of segments
$R[i,j]\in T(R[1,N])$. In particular, $R^{s}(1,N)=R[1,N]$. Set
$\mathbb{P}^{m}(i_R,j_R)$ and $\mathbb{P}^{s}(i_R,j_R)$ be the
probability of $R^{m}(i,j)$ and $R^{s}(i,j)$, respectively. Then we
have
$\mathbb{P}^{b}(i,j)=\mathbb{P}(L_1)+\mathbb{P}(L_2)+\mathbb{P}(L_3)
+\mathbb{P}(L_4)+\mathbb{P}(L_5)$, where
\begin{eqnarray*}
\mathbb{P}(L_{1})&=&\sum_{h,\ell}\mathbb{P}^{s}(h,\ell)\frac{Q^{s}(h,i-1)
Q^{b}(i,j)}{Q^{s}(h,\ell)}\\
\mathbb{P}(L_{2})&=&\sum_{h,\ell}\mathbb{P}^{b}(h,\ell)\frac{Q^{b}(i,j)
e^{-{\sf Int}(i,j;h,\ell)/kT}}{Q^{b}(h,\ell)}\\
\mathbb{P}(L_{3})&=&\sum_{h,\ell}\mathbb{P}^{b}(h,\ell)\frac{Q^{b}(i,j)
e^{-(\alpha_1+2\alpha_2+(\ell-j-1)\alpha_3)/kT}Q^{m}(h+1,i-1)}
{Q^{b}(h,\ell)}\\
\mathbb{P}(L_4)&=&\sum_{h,\ell}\mathbb{P}^{m}(h,\ell)\frac{Q^{b}(i,j)
e^{-((i-h+\ell-j)\alpha_3+\alpha_2)/kT}}{Q^{m}(h,\ell)}\\
\mathbb{P}(L_5)&=&\sum_{h,\ell}\mathbb{P}^{m}(h,\ell)\frac{Q^{b}(i,j)
e^{-((\ell-j)\alpha_3+\alpha_2)/kT}Q^{m}(h,i-1)}{Q^{m}(h,\ell)}.
\end{eqnarray*}
Accordingly, the recurrence formulae for $\mathbb{P}^{m}(i,j)$ and
$\mathbb{P}^{s}(i,j)$ are given as follows:\\
\begin{eqnarray*}
\mathbb{P}^{m}(i,j)&=&\mathbb{P}(M_1)+\mathbb{P}(M_2)\\
\mathbb{P}(M_1)&=&\sum_{h,\ell}\mathbb{P}^{b}(i-1,\ell)
\frac{e^{-(\alpha_1+2\alpha_2+(\ell-1-h)\alpha_3)/kT}
Q^{b}(j+1,h)Q^{m}(i,j)}{Q^{b}(i-1,\ell)}\\
\mathbb{P}(M_2)&=&\sum_{h,\ell}\mathbb{P}^{m}(i,\ell)
\frac{Q^{m}(i,j)Q^{b}(j+1,h)}{Q^{m}(i,\ell)}\\
\mathbb{P}^{s}(i,j)&=&\sum_{h,\ell}\mathbb{P}^{s}(i,\ell)
\frac{Q^{s}(i,j)Q^{b}(j+1,h)}{Q^{s}(i,\ell)}.
\end{eqnarray*}

\subsection{Base pairing probabilities for joint structures}
\begin{figure}[ht]
\centerline{%
\epsfig{file=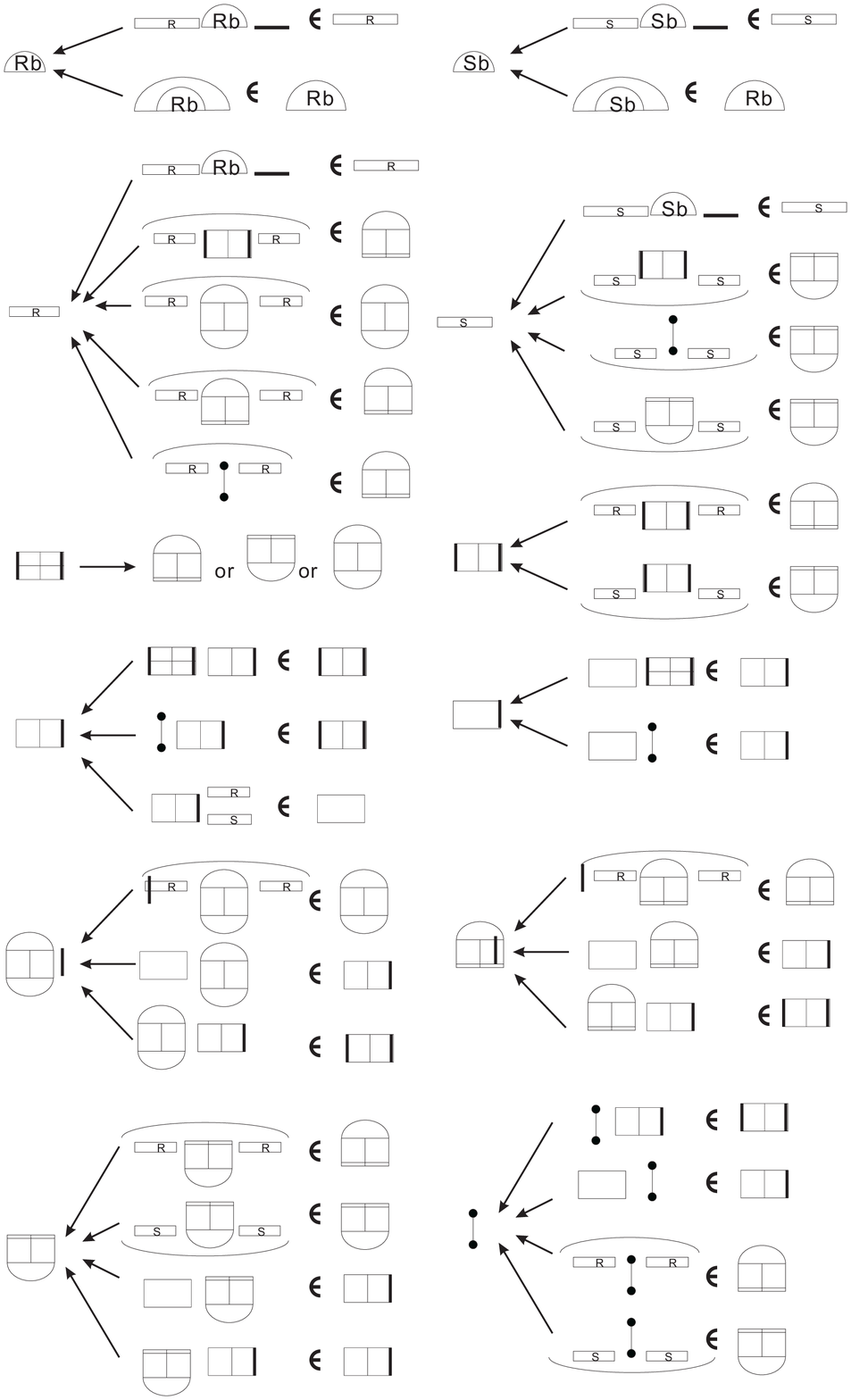,width=0.7 \textwidth}\hskip15pt
 }
\caption{\small Illustration of the ``inverse" grammar, obtained by
back-tracing in the decomposition tree. }\label{F:inversegrammar}
\end{figure}

Following the basic strategy, we first express the BPP via the
probabilities of particular substructures. In the following, we
abbreviate $J(1,N;1,M)$ by $J$. In order to calculate
$\mathbb{P}(i_{R},j_{R})$, let $\Sigma_{1}=\{J\mid (R[i],R[j])\in
J\}$, we consider the parent of $(R[i],R[j])$ in the $T(J)$ and
accordingly obtain
\begin{equation}
\Sigma_{1}=\{J\mid R^{b}[i,j]\in T(J)\}\cup \bigcup_{h,\ell}\{J\mid
J_{\bigtriangledown}(i,j;h,\ell)\in T(J)\}\cup \bigcup_{h,\ell}
\{J\mid J_{\square}(i,j;h,\ell)\in T(J)\},
\end{equation}
which immediately leads to
\begin{equation}\label{E:Rtt}
\mathbb{P}(i_{R},j_{R})=\mathbb{P}^{b}(i_{R},j_{R})+
\sum_{h,\ell}\mathbb{P}_{\bigtriangledown}
(i,j;h,\ell)+\sum_{h,\ell}\mathbb{P}_{\square}(i,j;h,\ell),
\end{equation}
where
\begin{eqnarray}
\mathbb{P}_{\bigtriangledown}(i,j;h,\ell)&=&\mathbb{P}^{E}_{\bigtriangledown}
(i,j;h,\ell)+\mathbb{P}^{M}_{\bigtriangledown}(i,j;h,\ell)+\mathbb{P}^{K}
_{\bigtriangledown}
(i,j;h,\ell)+\mathbb{P}^{F}_{\bigtriangledown}(i,j;h,\ell),\\
\mathbb{P}_{\square}(i,j;h,\ell)&=&\mathbb{P}^{E}_{\square}
(i,j;h,\ell)+\mathbb{P}^{M}_{\square}(i,j;h,\ell)+\mathbb{P}^{K}_{\square}
(i,j;h,\ell)+\mathbb{P}^{F}_{\square}(i,j;h,\ell).
\end{eqnarray}
Analogously, for $\mathbb{P}(i_{S},j_{S})$ we set
\begin{equation}
\Sigma_{2}=\{J\mid S^{b}[h,\ell]\in T(J)\}\cup \bigcup_{i,j}\{J\mid
J_{\bigtriangleup}(i,j;h,\ell)\in T(J)\},
\end{equation}
and obtain
\begin{equation}\label{E:Rtt1}
\mathbb{P}(i_{S},j_{S})=\mathbb{P}^{b}(h_{S},\ell_{S})+\sum_{h,\ell}\mathbb{P}_
{\bigtriangleup}(h,\ell;i,j),
\end{equation}
where
\begin{equation}
\mathbb{P} _{\bigtriangleup}
(h,\ell;i,j)=\mathbb{P}^{E}_{\bigtriangleup}(h,\ell;i,j)+\mathbb{P}^{M
}_{\bigtriangleup}(h,\ell;i,j)+\mathbb{P}^{K}
_{\bigtriangleup}(h,\ell;i,j)+\mathbb{P}^{F}
_{\bigtriangleup}(h,\ell;i,j).
\end{equation}
We remark that the expressions for the BPP $\mathbb{P}(i_{R},j_{R})$
and $\mathbb{P}(i_{S},j_{S})$ are not symmetric. This is due to the
fact that in our decomposition routines always the outer arcs
contained in $R$ are given preference. In other words, the asymmetry
is a result of our particular construction. Finally, we calculate
the binding probability of an exterior arc $(R[i],S[j])$. Since
$(R[i],S[j])$, being a tight structure of type $\circ$, is already
substructure, we can skip the first two steps of the basic strategy.
In order to compute the binding probabilities of both: interior and
exterior arcs, the key is to employ an ``inverse'' grammar induced
by tracing back in the decomposition tree as displayed in
Fig.~\ref{F:inversegrammar}. By virtue of this backtracking, we
obtain the recurrence formulae in analogy to the case of secondary
structures, discussed above.

\section{Synopsis}\label{S:syn}

In this paper we derive the partition function and the base pairing
probabilities of RNA interaction structures. Furthermore we present
the algorithm {\sf rip} that computes the partition function and the
base pairing probabilities in $O(N^4M^2)+O(N^2M^4)$ time and
$O(N^2M^2)$ space.

While the partition function is due to \cite{Backhofen} our construction
is independently derived and based on two ideas: the concept of tight
structure in Section~\ref{S:Introduction} and the decomposition tree,
presented in Section~\ref{S:grammar}.
We did however, adopt the notions of kissing and hybrid loops from
\cite{Backhofen}.
The derivation of the base pairing probabilities for joint structures is
new. Here the key idea is to express the latter via energy-wise
``quantifiable'' substructures, that are contained in the decomposition
tree. We discussed that in contrast to the computation of the base pairing
probabilities of secondary structures, the specific construction of
the unique grammar factors in. As a result, being a joint substructure
containing a certain base pair, is not the correct criterion any more.
Only those substructures that are obtained via tracing back in the
decomposition tree contribute to the base pairing probability.

The complete set of partition function recursions and all details on
the particular implementation of {\sf rip} can be found at
\begin{equation*}
{\tt http://www.combinatorics.cn/cbpc/rip.html}
\end{equation*}

Finally, we also compute the generating function of joint structures. The
analysis of this function is beyond the scope of this paper and can be
found as supplemental material at the above web-site.

\bibliography{rip}
\bibliographystyle{plain}
\end{document}